\documentclass[a4paper,12pt]{amsart}

\usepackage{fullpage}
\usepackage[driverfallback=dvipdfm]{hyperref}
\usepackage[all]{xy}
\usepackage{tikz}
\usepackage{amsaddr}
\usepackage{amsthm}
\usepackage{color}
\usepackage{enumerate}

\newtheorem{thm}{Theorem}[section]
\newtheorem*{thm*}{Theorem}
\newtheorem{cor}[thm]{Corollary}
\newtheorem{lem}[thm]{Lemma}
\newtheorem{prop}[thm]{Proposition}

\theoremstyle{definition}
\newtheorem{defn}[thm]{Definition}

\newtheorem{nt}[thm]{Notation}
\newtheorem{rem}[thm]{Remark}

\newtheorem*{rem*}{Remark}

\newtheorem*{theoremaux}{Theorem \theoremauxnum}
\gdef\theoremauxnum{1}

\def \p{{\mathbb P}}          

\def\N{{\mathbb N}}           
\def\R{{\mathbb R}}           
\def\O{{\mathcal O}}          
\def\Og{{\mathcal O}_{\text{gr}}}          
\def\G{{\mathcal G}}          
\def\e{\varepsilon}           
\def\f{\varphi}           

\def\Aut{\operatorname{Aut}}
\def\Out{\operatorname{Out}}
\def\Lip{\operatorname{Lip}}
\def\Hor{\operatorname{Hor}}
\def\PL{\operatorname{Str}}
\def\opt{\operatorname{opt}}
\def\wopt{\operatorname{weakopt}}
\def\vol{\operatorname{vol}}
\def\Min{\operatorname{Min}}

\def\TT{\operatorname{TT}}

\def\rank{\operatorname{rank}}
\def\core{\operatorname{core}}

\def\wt{\widetilde}
\def\ul{\underline}

\title[Displacements of automorphisms II]{
Displacements of automorphisms of free groups II:
Connectivity of level sets and decision problems}
\author{Stefano Francaviglia}
\address{Dipartimento di Matematica of the University of
Bologna}
\email{stefano.francaviglia@unibo.it}
\author{Armando Martino}
\address{Mathematical Sciences, University of Southampton }
\email{A.Martino@soton.ac.uk}
\begin{document}


\begin{abstract}
  This is the second of two papers in which we investigate  the properties  of displacement
  functions  of automorphisms of free groups (more generally, free products) on the
  Culler-Vogtmann Outer space $CV_n$ and its simplicial bordification. We develop a theory for
  both reducible and irreducible autormorphisms. As we reach the bordification of $CV_n$ we
  have to deal with general deformation spaces, for this reason we developed the theory in such
   generality. In first paper~\cite{FMpartI} we studied general properties of the
  displacement functions, such as well-orderability of the  spectrum and the topological
  characterization of min-points via partial train tracks (possibly at infinity). 

  This paper is devoted to proving that for any automorphism (reducible or not) any level set of
  the displacement function is connected. Here, by the ``level set" we intend to indicate the set of points displaced by \textit{at most  } some amount, rather than exactly some amount; this is sometimes called a ``sub-level set".

  As an application, this result provides a stopping procedure for
  brute force search algorithms in $CV_n$. We use this to reprove two known algorithmic results: the conjugacy problem for irreducible automorphisms and detecting irreducibility of automorphisms. 

\

Note: the two papers were originally packed together in the preprint~\cite{FMlevelset}
We decided to split that paper following the recommendations of a referee.

\end{abstract}
\maketitle
\tableofcontents

\section{Introduction}

We consider $F_n$ the free group of rank $n$, usually with a basis $B$ (a free generating
set). We are interested in the automorphism group, $\Aut(F_n)$ and the Outer automorphism
group, which is defined as $\Out(F_n)= \Aut(F_n)/\operatorname{Inn}(F_n)$. 

That said, as the 
reader will notice, in this paper all results are about general deformation spaces, and our statements are of the form ``let $[\phi]\in \Out(\Gamma)$'' or ``let
$X\in\overline{\Og(\Gamma)}^\infty$'' and so on. Let's briefly explain the notation and why we need to work in such generality. The reason is that classical Culler-Vogtmann space $CV_n$ is perfect for studying
irreducible automorphisms, but if one is interested in possibly reducible automorphisms, some
more general space is needed. If for instance an automorphism $\phi$ is represented by a
simplicial map $f$ on a finite graph $X$, it may happen that in $X$ we have a collection of
subgraphs $A_1,...,A_k$ so that $\cup_i A_i$ is preserved by $f$. In this case it may be
necessary to study both the invariant collection and the quotient obtained by collapsing any
$A_i$ to a point. So the typical object we have to deal with is a deformation space of finite
unions of graphs of groups. Concretely, our proofs boil down to induction proofs where the inductive step needs to deal with both the (disconnected) collection $\cup_i A_i$ and the map(s) that $f$ induces there, as well as the quotient graph of groups obtained by `collapsing' the $A_i$ in $X$, but keeping track of the fundamental group of the collapsed part; this leads to a graph of groups with trivial edge groups. So, even though our main focus is $CV_n$, it turns out to be no more complicated to deal with arbitrary free products and their deformation spaces, and our proofs need in fact to deal with the case of a finite graph of groups, with trivial edge groups, which {\em may not be connected}. This is what $\Gamma$ refers to. We direct the reader to Section~\ref{s3}, and in particular Remark~\ref{rem:gamma} for more discussion on this.

Nevertheless, since our general theorems specialise to results about classical $CV_n$ and $\Out(F_n)$,  
in this introduction we will stick as much as possible to  that classical setting.

\medskip

In recent years there has been a great deal of attention given to the Lipschitz metric on $CV_n$, see \cite{MR2862155}, \cite{MR2863547}, \cite{BestvinaBers} for instance. It has been considered even more generally in \cite{MR3342683}.

In the first part, \cite{FMpartI}, we proved results concerning the Lipschitz metric on a class of deformation spaces, of which a key example is the Culler-Vogtmann space of a free group, $CV_n$. We showed that, given an automorphism of a free group, the points of minimal displacement - for a given automorphism, the distance between a point in $CV_n$ and its image - correspond to the points which support partial train track maps, thus generalizing known results about irreducible automorphisms. 

In \cite{MR1396778} it is shown that, in the irreducible case, these points of minimal displacement (equivalently, the points which support train track maps) form a connected subset of $CV_n$ and this is used to solve the conjugacy problem. Our results here arise out of a desire to generalize those results to the reducible case, and we also employ Peak Reduction as a key tool.   

The generalization of this result for arbitrary, possible reducible, automorphisms, requires
some care, however. To start with, given an automorphism $\phi$, one can define the infimum
over all displacements of points in $CV_n$, to obtain $\lambda(\phi)$. However, in general
there might exist no points in $CV_n$ which are displaced by this amount. Our point of view is
to pass to the simplicial bordification of $CV_n$, otherwise known as the free splitting
complex, $\mathcal{FS}_n$. One can define displacements for points in $\mathcal{FS}_n$, though
in some cases these will be infinite. (A point in $CV_n$ is a marked graph, and a point in
$\mathcal{FS}_n$ arises by collapsing a subgraph. These induced points will have finite
displacement exactly when the subgraphs are $\phi$-invariant\footnote{See~\cite{FMpartI} or  Section~\ref{s3.4} for more details on this point.}). However, the infimum of all displacements of points in $\mathcal{FS}_n$ will, in general, be less than those in $CV_n$.

Bearing in mind these complications, and the fact that in the whole paper we work with more general
deformation spaces, our main Theorem, which is a special case
of Theorem~\ref{tconnected}, is the following: 

\bigskip
\parbox{0.9\textwidth}{{\bf Theorem} (Connectivity of Level Sets). Let $[\phi] \in \Out(F_n)$. Let $\lambda(\phi)$ be the infimum of displacements, with respect to the Lipschitz metric, of all points in $CV_n$. Then the set of points of $\mathcal{FS}_n$ which are displaced by exactly $\lambda(\phi)$,
is connected. }

\bigskip

\begin{rem*}
	As stated in the abstract we generally intend the ``level set" to be the set of points displaced by \textit{at most} some amount; this is sometimes referred to as a ``sub-level set". The subsequent Theorem has precisely this kind of statement. Hence the statement above is more properly a statement about the minimally displaced set, although our proofs deal with both at the same time.  
	
	However, note that $\lambda(\phi)$ is the infimum of displacements in $CV_n$; however, it might not be the infimum of displacements of points in $\mathcal{FS}_n$. 
\end{rem*}

\bigskip

Moreover, our techniques allow us to {\em regenerate} paths from $\mathcal{FS}_n$ to $CV_n$ without disturbing the displacements by very much. Hence, as part of the same Theorem~\ref{tconnected}, we also prove: 

\bigskip

\parbox{0.9\textwidth}{{\bf Theorem} Let $[\phi] \in \Out(F_n)$. Let $\lambda(\phi)$ be the infimum of displacements, with respect to the Lipschitz metric, of all points in $CV_n$. Then, for any $\e >0$ the set of points of $CV_n$ which are displaced by at most $\lambda(\phi)+\e$,
	is connected. }

\bigskip

In the case where the automorphism is irreducible, there are points in $CV_n$ which are displaced by exactly the minimum, $\lambda(\phi)$. Moreover, every point on the boundary has infinite displacement (Remark~\ref{displacement at infinity}) and hence the connectivity of the level set becomes a statement about $CV_n$, as in Corollary~\ref{tconnectedirred}:  

\bigskip

\parbox{0.9\textwidth}{{\bf Corollary} Let $[\phi] \in \Out(F_n)$ be irreducible. Let $\lambda(\phi)$ be the infimum of displacements, with respect to the Lipschitz metric, of all points in $CV_n$. Then the set of points of $CV_n$ which are displaced by $\lambda(\phi)$,
	is connected. }

\bigskip

\begin{rem}
	Given an automorphism, $\phi$, of the free group, one can construct a relative train track representative for $\phi$. The quantity, $\lambda(\phi)$ is then simply the maximum Perron-Frobenius eigenvalue of any stratum. 
	
	More generally, if we are given a $\phi$-invariant free factor system, then one can build a relative train track representative of $\phi$ which sees this free factor system as an invariant subgraph. There is a corresponding deformation space where one collapses this subgraph, and the minimum displacement in that deformation space is the maximum  Perron-Frobenius eigenvalue of any stratum  {\em above} the invariant subgraph. 
	
	We can think of $\mathcal{FS}_n$ as a union of such deformation spaces, with the displacements being infinite when the collapsed object is not invariant. This is why the minimum displacement in $\mathcal{FS}_n$ need not be equal to that in $CV_n$ - they are different if one can collapse an invariant subgraph which carries all the maximum Perron-Frobenius eigenvalues.

	A simple example is the following. Consider this automorphism, $\phi$, of the free group on $a,b,c$:

	$$
	\begin{array}{rcl}
    c & \mapsto & ca\\	
	b & \mapsto & ba \\
	a & \mapsto & aba \\
	\end{array}
	$$
This is then a relative train track map, with two strata, the bottom one given by $a,b$ and the top one by $c$.

Let $\lambda$ be the larger eigenvalue of the matrix 
$$
\begin{bmatrix}
2 & 1 \\
1 & 1\\
\end{bmatrix}.
$$	
This is the Perron-Frobenius  eigenvalue of the bottom stratum, with the top stratum having 1 as its Perron-Frobenius  eigenvalue.  It is then easy to see that $\lambda(\phi)= \lambda$, but there are points in $\mathcal{FS}_n$  which are fixed by $\phi$ and so have multiplicative displacement $1$; namely, take the point obtained by collapsing $a,b$. That is, the graph of groups with one edge, one vertex, a trivial edge group and a vertex group generated by $a$ and $b$. 
\end{rem}

\medskip

Naturally, since our results generalize those of \cite{MR1396778}, we obtain a solution of the conjugacy problem for irreducible automorphisms in the same way. However, it seems that our techniques allow for a more elementary interpretation, and also opens up the possibilty for attempting the algorithm in the reducible case. However, there are further complications that arise in the reducible case, due to the fact that the minimally displaced set enters the thin part, and so we do not easily obtain bounds on the number of points we need to enumerate.  

\medskip

In any case, we can describe this algorithm in the irreducible case, with explicit constants, rather straightforwardly. Moreover, we also provide an algorithm to detect irreducibility; this result was first proved in \cite{MR3194747}
and improved in \cite{kapovichpolyred}  (also, see \cite{ctmaps} and \cite{MR2763779} which
give another algorithm for detecting irreducibility).

\medskip

Finally, it may be worth noting that, thank to the generality of Theorem~\ref{tconnected}, our
algorithms easily generalise to a class of groups bigger than just free groups; concretely,
groups of the kind $G_1*\dots*G_p*F_n$ where the $G_i$ are finite groups (see Theorem~\ref{finite} and Section~\ref{s2.1}).

\bigskip

\noindent
\textsc{Acknowledgements:} We would like to thank both the Universit\`{a} di Bologna and the
Universitat Polit\'{e}cnica de Catalunya, for their hospitality during several visits. We would
also like to thank the referee of the earlier paper~\cite{FMlevelset}, as well as the current referee for their patience and extraordinary efforts in improving this paper, with their many helpful comments and suggestions.

\section{Algorithms}

 In order to motivate the detailed discussion which follows, we provide here the two algorithms for solving conjugacy in the irreducible case and for detecting irreducibility. We present these algorithms as naively as possible, in order to make them more accessible. That is, one could understand and implement them without any knowledge of the Lipschitz metric, Culler-Vogtmann space or partial train track maps. As such we have made no attempt to streamline the algorithms in any way; they are brute force searches in an exponential space. 

However, we would stress that our point of view is fundamentally that these procedures would be better run as path searches in Culler-Vogtmann space, enumerating optimal maps and calculating displacements via candidates. That abundance of terminology would make the algorithms much harder to describe, so we instead translate everything to a more manageable setting; bases of $F_n$ and generating sets for $\Out(F_n)$. However, the technical point of view is more helpful in developing an intuition of the processes and is likely the way to vastly improve the algorithmic complexity. 

\medskip

Let us now describe our algorithms, whose correctness is proved at the end of the paper. First, we recall some terminology. In order to work algorithmically with $\Out(F_n)$ we need a generating set. The best known of these is the set of Nielsen generators, but it is more convenient for us to work with the following:

\begin{defn}[CMT Automorphisms, \cite{MR721773} and \cite{MR748994}]\label{defCMT}
A {\em CMT} automorphism of $F_n$ is one that is induced by a change of maximal tree. More
precisely, fix a basis, $B$, of $F_n$. Let $R=R_B$ be the marked rose corresponding to $B$; that is, $R$ is a graph with one vertex $v$ and $n$ edges called petals, and we have a fixed isomorphism between $F_n$ and $\pi_1(R,v)$ where each element of $B$ corresponds to a petal of $R$ under this isomorphism.

Let $X$ be a graph with fundamental group of rank $n$, and let $T, T'$ be two maximal trees of
$X$. Collapsing $T$ and $T'$ we obtain two roses $R_T$ and $R_{T'}$.  Let $\rho_T, \rho_{T'}$
be the corresponding projections from $X$ to $R_T,R_{T'}$, and let $\alpha_T,\alpha_{T'}$ be
homeomorphisms from $R$ to $R_T,R_{T'}$ respectively.  Then the (outer) automorphism induced by
changing the maximal tree from $T$ to $T'$ is the (homotopy class of the) map
$\alpha^{-1}_{T'}\rho_{T'} {\rho_T}^{-1}\alpha_T:R\to R$, where the inverse denotes a homotopy inverse.

Thus the set of CMT automorphisms of $F_n$, relative to $B$, is the set of all such change of maximal tree automorphisms.  
\end{defn}

The set of CMT automorphisms relative to $B$ includes all Whitehead automorphisms, (see \cite[Theorem~5.5]{MR721773}, and \cite{MR1922276}) and is a finite set which generates $\Out(F_n)$. 
Note also that in case $T=T'$, by varying $\alpha_T,\alpha_{T'}$ we obtain all graph automorphisms of $R$,
including inversions of generators, which therefore are CMT automorphisms. Also, note that the
property of being CMT, depends on a fixed chosen basis, $B$, of $F_n$ (the petals of $R$).

\begin{rem}
	The definition of CMT automorphisms just given is close to that given in \cite{MR1922276}, but there is an alternative definition via $CV_n$ as follows. We call two marked roses, $R_1, R_2$ adjacent if there is a simplex, $\Delta$, in $CV_n$, admitting faces, $\Delta_1$ and $\Delta_2$ such that $R_i$ belongs to $\Delta_i$. This is equivalent to saying there is a marked graph $X$, admitting two maximal trees, $T_1, T_2$ whose collapse produces the marked graphs, $R_1, R_2$, respectively. Then, 
	
	$$
	CMT_R(F_n) = \{ [\phi] \in \Out(F_n) \ :  \ \phi(R) \text{ is adjacent to } R \}.
	$$ 
	
This is the same as the previous definition by setting $R=R_B$, as above. 	
\end{rem}

\medskip

Next we need a notion of size of an automorphism, which will provide a termination criterion for our algorithms. 

\begin{defn}
Let $[\phi] \in \Out(F_n)$, and let $B$ be a basis of $F_n$. Define $|| \phi ||_B$ to be $\sup_{1
  \neq g \in F_n} \frac{|| \phi g||_B}{||g||_B}$, where $||g||_B$ denotes the cyclic reduced length of
$g$ with respect to $B$. This supremum is a maximum and is realised by an element of cyclic
length $\leq 2$ (see Lemma~\ref{sausagelemma}).  
\end{defn}

\begin{rem}\label{remfin}
Note that for any constant, $C$, there are only finitely many $[\phi] \in \Out(F_n)$ such that $||\phi||_B \leq C$ (see \cite{Francaviglia2009}, Lemma 4.10). 

This also follows since $||\phi||_B$ is really $\Lambda(R, \phi(R))$ in disguise (see section~\ref{s3.4}) where $R=R_B$ is the uniform marked rose corresponding to $B$ with all edges length $1/n$ (so that $R$ has volume 1). The remark then follows from the fact that, for any given $C$, there are only finitely many (marked, volume 1) roses, $R_1$ such that $\Lambda(R, R_1) \leq C$, and the stabiliser of any point, and in particular $R$, is finite.

%
%
\end{rem}

Our first application is then as follows. (See Section~\ref{appl} for the proof.)

\begin{thm} \label{conjirred}
The following is an algorithm to determine whether two irreducible automorphisms are conjugate. 

\medskip

\noindent
Let $[\phi], [\psi]$ be two irreducible outer automorphisms of $F_n$, and $B$ a basis of $F_n$.
\begin{itemize}
\item Choose any $\mu > \max \{ ||\phi||_B, ||\psi||_B \}$. 
\item Inductively construct a finite set, $S=S_{\phi, \mu}$, as follows (which depends on both $\phi$ and $\mu$): 
\begin{itemize}
\item Start with $S_0 = \{ \phi \}$. 
\item Set $K= n(3n-3) \mu^{3n-1}$. 
\item Inductively put $S_{i+1}$ to be all possible automorphisms $ \zeta \phi_i \zeta^{-1}$,
  where $\phi_i$ is any element of $S_i$, $\zeta$ is any CMT automorphism, subject to the
  constraint that $||   \zeta \phi_i \zeta^{-1} ||_B \leq K$. (Since the identity is a CMT
  automorphism according Definition~\ref{defCMT}, we have $S_{i-1} \subseteq S_i$). 
\item End this process when $S_i = S_{i+1}$, and let this final set be $S$. 
\end{itemize} 
\item Then $\psi$ is conjugate to $\phi$ if and only if $\psi \in S$. 
\end{itemize}
\end{thm}

\medskip

Of course, one would like to also be able to decide when an automorphism is irreducible when it is given by images of a basis, for instance. In order to do so, we recall the definition of irreducibility. 

\begin{defn}[See \cite{BestvinaHandel}]
An (outer) automorphism, $[\phi]$ of $F_n$ is called {\em reducible} if there are free factors, $F_{n_1}, \ldots, F_{n_k}, F_{n_\infty}$ such that $F_n = F_{n_1} * \ldots F_{n_k} * F_{n_\infty}$ and each $\phi(F_{n_i})$ is conjugate to $F_{n_{i+1}}$ (subscripts taken modulo $k$). If $k=1$ we further require that $F_{n_\infty} \neq 1$. (In general $\phi(F_{n_\infty})$ is not conjugate to $F_{n_\infty}$). Otherwise $[\phi]$ is called irreducible. 

Equivalently, $[\phi]$ is reducible if it is represented by a homotopy equivalence $f$, on a
core graph $X$, such that $X$ has a proper subgraph $X_0$, with non-trivial fundamental group, such that $f(X_0) = X_0$. (Being represented by $f$ means that there is an isomorphism, $\tau: F_n \to \pi_1(X)$ such that $\phi = \tau^{-1} f_* \tau$).  
\end{defn}

We add the following, which constitutes an obvious way that one can detect irreducibility by inspection.  

\begin{defn}
Consider $F_n$ with basis $B$ and let $[\phi]$ be an outer automorphism of $F_n$. We say that $[\phi]$ is visibly reducible with respect to $B$, or simply visibly reducible, if there exist disjoint subsets $B_1, \ldots, B_k$ of $B$ such that $\phi(\langle B_i \rangle)$ is conjugate to $\langle B_{i+1} \rangle$ (with subscripts taken modulo $k$). If $k=1$ we also require that $B_1 \neq B$. 

More generally, we say that a homotopy equivalence on the rose is visibly reducible if it is visibly reducible with respect to the basis given by petals. 
\end{defn}

This is, in fact, easy to check by classical methods due to Stallings, \cite{Sta}.

\begin{lem}
If $[\phi]$ is visibly reducible, then it is reducible. Moreover, there is an algorithm to determine if $[\phi]$ is visibly reducible with respect to $B$. 
\end{lem}
\proof
The first statement is clear, since each subset of a basis generates a free factor, and disjoint subsets generate complementary free factors. Since there are only finitely many subsets to check, we simply need to determine if the conditions that $\phi(\langle B_i \rangle)$ is conjugate to $\langle B_{i+1} \rangle$ hold. But this can readily be checked since two subgroups of a free group are conjugate if and only if the core of their Stallings graphs are equal, \cite{Sta}. 
\qed

\medskip
We can now describe our second algorithm. (See Section~\ref{appl} for the proof.)

\begin{thm} \label{detectirred}
The following is an algorithm to determine whether or not an outer automorphism of $F_n$ is irreducible.

\medskip

\noindent
Let $[\phi]\in\Out(F_n)$, and $B$ a basis of $F_n$. Construct $S=S_{\phi}$ as above. Namely, 
\begin{itemize}
\item Choose any $\mu > ||\phi||_B$. 
\item Inductively construct the finite set, $S=S_{\phi, \mu}$: 
\begin{itemize}
\item Start with $S_0 = \{ \phi \}$. 
\item Set $K= n(3n-3) \mu^{3n-1}$. 
\item Inductively put $S_{i+1}$ to be all possible automorphisms $ \zeta \phi_i \zeta^{-1}$,
  where $\phi_i$ is any element of $S_i$, $\zeta$ is any CMT automorphism, subject to the
  constraint that $||   \zeta \phi_i \zeta^{-1} ||_B \leq K$.
  (Since the identity is a CMT automorphism according Definition~\ref{defCMT}, we have  $S_{i-1}
  \subseteq S_i$).  
\item End this process when $S_i = S_{i+1}$, and let this final set be $S$. 
\end{itemize} 
\item Let $S^+$ be the set of all possible automorphisms $ \zeta \phi_i \zeta^{-1}$, where $\phi_i$ is any element of $S$, $\zeta$ is any CMT automorphism, with no other constraint. 
\item If some $\psi \in S^+$ is visibly reducible with respect to $B$, then $\phi$ is reducible. Otherwise, $\phi$ is irreducible. 
\end{itemize}
\end{thm}
\begin{rem}\label{remS}
  In both Theorems~\ref{conjirred}, and~\ref{detectirred}, the set $S$ is a subset of the
  set of automorphism classes with $||\phi||_B\leq K$, which is finite
  (Remark~\ref{remfin}).
  Therefore both algorithms stop in a finite, effectively computable, time. 
\end{rem}

We also explain, in Section~\ref{s2.1}, how to implement essentially the same algorithms in the case where one has a free product of finite groups with a free group. 

\begin{thm}
	\label{finite}
	Let $G = G_1* \ldots *G_p*F_n$ be a free product where the $G_i$ are finite groups and $F_n$ is a free group of rank $n$.  Let $\G = \{ \{ G_i\}, n\}$ be the free splitting induced from the finite groups $G_i$. Then the following problems are algorithmically decidable: 
	\begin{itemize}
		\item Deciding whether a given $[\phi] \in \Out(G)$ is irreducible (relative to $\G$),
		\item Deciding whether two $\G$-irreducible automorphisms, $[\phi],[\psi]$ are conjugate in $\Out(G)$.  
	\end{itemize}
\end{thm}

\begin{rem}
	Note that any automorphism of $G$ preserves $\G$, so $\Out(G) = \Out(\G)$ in this case. 
\end{rem}

\newpage

\section[Preliminaries]{Preliminaries and notation (from~\cite{FMpartI})}\label{s3}
Throughout the paper, we use the definitions and notation of~\cite{FMpartI}.
We briefly recall them here, referring the reader to~\cite{FMpartI} for a detailed discussion.

Before of that, we wish to recall the reasons for giving new definitions and working in a so
general setting. Our principal motivation was to study outer automorphisms of free groups that are
possibly reducible. This naturally leads to consider simplicial bordifications of 
Culler-Vogtmann Outer spaces. 
Namely, if $\Gamma$ is a marked graph with fundamental group $F_n$ -- the rank-$n$
free group -- then any automorphism $\phi:F_n\to F_n$ can be represented by a simplicial map
$f:\Gamma\to \Gamma$. When $\phi$ is reducible, it may happen that $\Gamma$ exhibits a
collection $\Gamma_1,\dots,\Gamma_k$ of $f$-invariant subgraphs. In order to study the
properties of $\phi$ it may help to collapse such sub-graphs to points. So one is naturally
lead to study two kind of deformation spaces: that of actions on trees with possibly non-trivial
vertex stabilisers (when we collapse the $\Gamma_i$'s) and product of such spaces (when we
consider the restriction of $f$ to the invariant collection).

Summing up, the typical object we need to understand is a disjoint union of metric trees, where
a group $G$ acts with possibly non-trivial vertex stabilisers. We therefore work in such
a general setting, as developed in~\cite{FMpartI},  but the reader is invited to keep in mind the case of $CV_n$ and its
bordification.

\subsection{Splittings, $\G$-trees, outer spaces, and automorphisms}

$F_n$ will always denote the free group of rank $n$. We will consider groups $G$ equipped with
free a splitting $G=G_1*\dots*G_p*F_n$. We do not assume $G_i$ is indecomposable, and our main
interest is indeed when $G$ itself is a free group.

\begin{defn}\label{d31r}
  Given a group $G$, a {\em free splitting} $\G$ of $G$ is a pair $(\{G_1,\dots,G_p\},n)$ where
  $\{G_i\}$ is a collection of subgroups  of $G$ (and $n\in \N$) such that $G=G_1*\dots*G_p*F_n$.
Two splittings $\G=(\{G_i\},n)$ and $\mathcal H=(\{H_i\},m)$ of $G$ are {\em of the same type} if $m=n$ and, up to
reordering and conjugacy of the $G_i$, they have the same factor subgroups. That is, we do not require the named (conjugacy class of the) free group factor at the end to be preserved. 
  The Kurosh {\em rank} of the splitting is $n+p$. We say that $\mathcal H$ is a {\em sub-splitting} of $\G$ if every $H_i$ decomposes as
$H_i=G_{j_1}*\dots*G_{j_{l_i}}*F_{s_i}$ and $n=m+\sum s_i$.
\end{defn}

\begin{rem}
  We admit the trivial splitting $G=F_n$, $(\emptyset, n)$. That is the splitting with no free factors groups. In this case our discussion will amount to considering the free group $F_n$ and the classical Culler-Vogtmann Outer space $CV_n$.
\end{rem}

\begin{rem}
  Free splittings are also referred to as {\em free factor systems} in the literature
  (originally introduced in \cite{BFH2}, and also used in \cite{handelmosher}, \cite{HM2} and
  \cite{horbez}). The viewpoint of \cite{FMpartI} and the present paper is that of taking a
  fixed free factor system - a free splitting - and  studying its deformation space. We refer
  to~\cite{FMpartI} for more details. We just notice here that a ``splitting'' in general refers to any action on a tree and the induced graph of groups decomposition, but no confusion should arise since all of the splittings we consider are ``free'', in the sense that the edge stabilisers in the tree are trivial (equivalently, the splitting which arises is a free factor system).

\end{rem}

\begin{defn}
  Given a group $G$ endowed with a free splitting $\G=(\{G_i\},n)$, a simplicial
  {\em $\G$-tree} is a metric simplicial tree, endowed with a faithful simplicial $G$-action via isometries, trivial
edge-stabilisers, and such that for every $G_i$ there is exactly one orbit of vertices whose
stabiliser is conjugate to $G_i$. Such vertices are called {\em non-free}. Other vertices
(those with trivial stabilisers) are called {\em free} vertices. 

A {\em $\G$-graph} is a finite connected metric graph of groups $X$ whose {\em topological}
fundamental group is $F_n$,  with trivial edge-groups, and endowed with a $G$-marking, that is,
there is a fixed isomorphism between its fundamental group -- {\em as graph of groups} -- and $G$, such
that the splitting given by vertex groups is equivalent to $\G$.

  If the splitting $\G$ of a group $G$ is clear from the context, we may use notation $G$-tree 
  instead of $\G$-tree. Same for graph.

  The rank of a $\G$-tree (or graph) is the Kurosh rank of the splitting (as defined in Definition~\ref{d31r}).

\end{defn}

A $\G$-tree is {\em mimimal} if it has no proper $G$-invariant sub-tree (in particular, it has
no free leaves, and $G$ acts without global fixed point). A graph of groups with trivial
edge-groups is a {\em core graph} if its leaves (if any) have non-trivial vertex group.  
Given a graph of groups $X$, with trivial edge groups and non-trivial fundamental group (as
graph of groups), we define $\core(X)$ to be the maximal core
  sub-graph of $X$. If $X$ has trivial fundamental group (as graph of groups) we define $\core(X)$ to be just a
  point of $X$. We say that $\core(X)$ is trivial when it is a point, namely when $X$ is 
  topologically a tree with at most one non-free vertex.

Bass-Serre theory provides a correspondence between minimal $\G$-tree and core $\G$-graphs, so
one can equivalently works either with trees or graphs. The equivalence tree-graph is made
explicit as follows: Given a minimal $\G$-tree, its quotient by the $G$-action is a core
$\G$-graph.

Two simplicial $\G$-trees are considered {\em equivalent} if there is a $G$-equivariant
isometry between them, and the corresponding notion of equivalence is given for graphs.

In some setting it will be more convenient
using trees, in others, graphs. For this purpose we introduce the following notation.
\begin{nt}[Tilde-underline notation]\label{tildeunderbar}
 Let $\G$ be a free splitting of a group $G$.
If $X$ is a $\G$-graph, then $\wt X$  denotes its universal covering, which is a $\G$-tree. As
usual, if $x\in X$ then $\wt x$ will denote  a lift of $x$ in $\wt X$. The same for subsets: if
$A\subset X$ is connected then $\wt A\subset \wt X$ is a connected component of the preimage of
$A$. On the converse situation, if $T$ is a $\G$-tree with finite edge-orbits, we denote by $\ul
  T$ the quotient $\G$-graph. Same notation for points and subsets.  So, $\widetilde{\ul X}=X$
  for both graphs and trees. 
\end{nt}

\begin{nt}
 If $X$ is a connected graph of
  groups with trivial edge groups, by a $X$-graph (or tree) we mean a $\pi_1(X)$-graph
  (resp. tree), that is to say, a   $\G$-graph (resp. tree) where $G$ is the fundamental group 
  of $X$ as graph of groups, and $\G$ is the  splitting of $G$ given by vertex groups.  

 If $\Gamma=\sqcup \Gamma_i$ is a
  disjoint finite union of finite graphs of groups 
  with trivial edge-groups, a $\Gamma$-graph is a disjoint finite union $X=\sqcup X_i$ of
  $\Gamma_i$-graphs (and a $\Gamma$-forest is a union of $\Gamma_i$-trees).  
\end{nt}

We introduce now the outer space of a splitting (see~\cite{FM13,GuirardelLevitt,FMpartI} for details).
  Let $G$ be a group and $\G$ be a splitting of $G$. The (projectivized) outer space of
  $G$, relative to the splitting $\G$, consists of (projective)
  classes of minimal simplicial metric $\G$-trees $X$ with no redundant vertex
  (i.e. free and two-valent) and such that the $G$-action is by isometries.\footnote{If $\G$ is the trivial splitting $G=F_n$, then $\O(\G)=CV_n$.} 

We use the notation $\O(G;\G)$ or simply $\O(\G)$ to indicate the outer space
of $G$ relative to $\G$. We use $\p\O(G;\G)$ (or simply  $\p\O(\G)$) to indicate the 
projectivized outer space.  For $X\in\O(\G)$ we define its volume $\vol(X)$ as the sum of
lengths of edges in $G\backslash X$. This is often referred to also as co-volume.
The volume-one slice of $\O(\G)$ is indicated by $\O_1(\G)$\footnote{We stress that the
  distinction between $\O(\G)$ and $\p\O(\G)$ is not crucial in our setting as we will mainly work with scale-invariant functions.}.

We defined $\O(\G)$ as a space of {\bf trees}, but we it will be often convenient to use {\bf graphs} $X$
so that $\wt X\in \O(\G)$. Clearly the two viewpoints are equivalent. We introduce the
following convention: when we want to consider spaces of graphs we add a ``lower $gr$'' to our
notation: $$\Og(\G)=\{\G\text{-graph } X: \wt X\in\O(\G)\}$$
The {\bf spaces $\O(\G)$ and $\Og(\G)$ are naturally identified via $X\leftrightarrow \wt X$.}
In particular, they are completely interchangeable in all statements.

If $X$ is a finite connected graph of groups with trivial edge-groups, and $\mathcal S$ is the
splitting of $\pi_1(X)$ given by vertex-groups, then we set 
$$\O(X)=\O(\pi_1(X);\mathcal S).$$

Let now $\G$ be a splitting of a group $G$, 
$X$ be a $\G$-graph, and $\Gamma=\sqcup_i\Gamma_i$ be a sub-graph of $X$ whose connected
components $\Gamma_i$ have non-trivial fundamental groups (as graphs of groups).
Then $\Gamma$ induces a sub-splitting $\mathcal S$ of $\G$ where the factor-groups $H_j$ are either
\begin{itemize}
\item the fundamental groups $\pi_1(\Gamma_i)$, or
\item the vertex-groups of non-free vertices in $X\setminus \Gamma$.
\end{itemize}
 In this case will use the notation
 $$\O(X/\Gamma):=\O(G;\mathcal S)\qquad \O(\Gamma):=\Pi_i\O(\Gamma_i)$$
(and similarly for $\Og$).
 We tacitly identify $X=(X_1,\dots,X_k)\in \O(\Gamma)$ with the labelled disjoint union
$X=\sqcup_i X_i$. So an element of $\O(\Gamma)$ can be interpreted as a metric $\Gamma$-forest.
The quotient of $\O(\Gamma)$ by the natural action of $\mathbb R^+$ is the projective outer
space of $\Gamma$, and it is denoted by $\p\O(\Gamma)$. (Thus $\p\O(\Gamma)$ is not the product
of the $\p\O(\Gamma_i)$'s.) The notion of volume extends to $\Gamma$-trees: If
$X=(X_1,\dots,X_k)\in\O(\Gamma)$ we set $\vol(X)=\sum_i\vol(X_i)$, and $\O_1(\Gamma)$ denotes
the volume-one slice of $\O(\Gamma)$. We extend our notation and define define $\O(X/A)$
and $\O(A)$ also to the case where $X$ is a non connected $\Gamma$-graph and $A\subset X$ is a sub-graph whose components have non-trivial fundamental groups.

\begin{nt}\label{not:gamma}
In what follows we use the following convention:
\begin{itemize}
\item $G$ will always be a group with  a splitting $\G=(\{G_1,\dots,G_p\},F_n)$;
\item  $\Gamma=\sqcup \Gamma_i$ will always mean that $\Gamma$ is a finite disjoint
union of finite graphs of groups $\Gamma_i$, each with trivial edge-groups and non-trivial fundamental group $H_i=\pi_i(\Gamma_i)$, each
$H_i$ being equipped with the splitting given by the vertex-groups.
\end{itemize}
\end{nt}

We set $$\rank(\Gamma)=\sum_i\rank(\Gamma_i).$$

\begin{rem}\label{rem:gamma}
One should think $G$-statements as referring to classical outer space $CV_n$, $\G$-statements
as referring to its simplicial bordification and deformation spaces of free products, and $\Gamma$-statements as general statements
about more general deformation spaces, that come into play along the way of our rank-inductive
strategy. More precisely, any $\Gamma$-statement specialises to a $\G$-statement (in the case where $\Gamma$ is connected),
to a $G$-statement ($\Gamma$ connected and trivial splitting), and to a $CV_n$-statement
($\Gamma$ is connected, the splitting is trivial, and $G=F_n$).

For this reason, the paper will contain mainly $\Gamma$-statements.
\end{rem}

\begin{nt}
\label{marked}

We will also consider moduli spaces with marked points. The moduli space of $\G$-trees with $k$ marked points $p_1,\dots,p_k$ (not necessarily distinct) is denoted by 
  $\O(G;\G,k)$ or simply $\O(\G,k)$. 
If $\Gamma=\sqcup_{i=1}^s \Gamma_i$, given $k_1,\dots,k_s\in\N$ we
  set $$\O(\Gamma,k_1,\dots,k_s)=\Pi_i\O(\Gamma_i,k_i).$$
\end{nt}

\medskip
We introduce now the group $\Aut(\Gamma)$.
The group of automorphisms of $G$  that preserve the set of conjugacy classes of the $G_i$'s is
  denoted by $\operatorname{Aut}(G;\G)$. We set
  $\operatorname{Out}(G;\G)=\operatorname{Aut}(G;\G)/\operatorname{Inn}(G)$

The group $\Aut(G;\G)$ acts on $\O(\G)$ by
changing the marking (i.e. the action), and $\operatorname{Inn}(G)$ acts trivially. Hence
$\operatorname{Out}(G;\G)$ acts on
$\O(\G)$. If $X\in\O(\G)$ and
$[\phi]\in\operatorname{Out}(G;\G)$ then $\phi X$ is the same metric tree as $X$, but the action
is $(g,x)\to \phi(g)x$. The action is simplicial and continuous with respect to both simplicial and
equivariant Gromov topologies. (See Section~\ref{s2} for details on simplicial structures). We
remark that despite the left notation, this is a right-action.

\medskip
We now extend the definition of $\Aut(G;\G)$ to the case of $\Gamma=\sqcup_i\Gamma_i$. We
denote by
$\mathfrak S_k$ the group of permutations of $k$ elements.

  Let $G$ and $H$ be two isomorphic groups endowed with splitting $\G:G=G_1*\dots G_p*F_n$ and
$\mathcal H:H=H_1*\dots H_p*F_n$. The set of isomorphisms from $G$ to $H$ that maps each $G_i$
to a conjugate of one of the $H_i$'s is denoted by $\operatorname{Isom}(G,H;\G,\mathcal H)$. If
splittings  are clear from the context we write simply
$\operatorname{Isom}(G,H)$. 

\begin{defn} \label{autgamma}
 For $\Gamma=\sqcup_{i=1}^k\Gamma_i$  as in Notation~\ref{not:gamma}, we set $$\Aut(\Gamma)=\{\phi=(\sigma,\phi_1,\dots,\phi_k):\ \sigma\in\mathfrak S_k\text{
and } \phi_i\in\operatorname{Isom}(H_i,H_{\sigma_i})\}.$$
$$\operatorname{Inn}(\Gamma)=\{(\sigma,\phi_1,\dots,\phi_k)\in\Aut(\Gamma): \sigma=id,
\phi_i\in\operatorname{Inn}(H_i)\}$$
$$\Out(\Gamma)=\Aut(\Gamma)/\operatorname{Inn}(\Gamma).$$
\end{defn}

The composition of $\Aut(\Gamma)$ is component by component defined as follows. Given
$\phi=(\sigma,\phi_1,\dots,\phi_k)$ and $\psi=(\tau,\psi_1,\dots,\psi_k)$ we have
$$\psi\phi=(\tau\sigma,\psi_{\sigma(1)}\phi_1,\dots,\psi_{\sigma(k)}\phi_k).$$
The group $\Aut(\Gamma)$ acts on $\O(\Gamma)$ in the natural way with kernel
$\operatorname{Inn}(\Gamma)$, so $\Out(\Gamma)$ acts on $\O(\Gamma)$. More precisely, if
$X=(X_1,...,X_k)\in\O(\Gamma)$ then $X_{\sigma(i)}$ becomes an $H_i$-tree --- denoted
$\phi_iX_{\sigma(i)}$ ---  via the pre-composition of $\phi_i:H_i\to H_{\sigma(i)}$ with the
$H_{\sigma(i)}$-action. We set $\phi X=(\phi_1X_{\sigma(1)},\dots,\phi_kX_{\sigma(k)})$. 
(See~\cite[Section 2]{FMpartI} for more details).

\subsection{Simplicial structure of outer spaces and its bordification (\cite[Sections~2.5 and~2.6]{FMpartI})}\label{s2}
The simplicial structure we are going to use is the usual one, that is, (open) simplices are
defined as follows: for any $X\in\O(\G)$, the set of $\G$-trees obtained from $X$ by varying
edge-orbit-lengths in  $(0,\infty)$, is an open simplex of $\O(\G)$, that we refer to as the open
simplex of $X$, and denote by $\Delta_X$. We notice that $\Delta_X\cap \O_1(\G)$ is the
standard open simplex in $\R^{\textrm{number of edge-orbits}}$, while $\Delta_X$ is its
positive cone (which is topologically still an open simplex, just one dimension bigger). On any open simplex
we  put the {\bf  Euclidean} sup-distance $d_\Delta^{Euclid}(X,Y)$ ($d_\Delta(X,Y)$ or $d(X,Y)$
for  short) 
$$d_\Delta^{Euclid}(X,Y)=d_\Delta(X,Y)=\max_{e\text{ edge}}|L_X(e)-L_Y(e)|.$$ 
Such definitions naturally extend to the case of $\Gamma=\sqcup_i\Gamma_i$.
(Note, however, that the simplicial structure of $\p\O(\Gamma)$ is not the product of the structures of
$\p\O(\pi_1(\Gamma_i))$.)

\begin{rem}
  The identification of a simplex $\Delta$ with a subset of $\R^m$, induces the notion
  of linear combination $sX+tY$ for any $X,Y\in\Delta$ and $s,t\geq0$. In particular, the
  convex combination $tX+(1-t)Y$ is well defined for any $t\in[0,1]$. We refer to the set
  $\overline{XY}:=\{tX+(1-t)Y, t\in[0,1]\}$ as the Euclidean segment between $X$ and $Y$.
\end{rem}

{\em Simplicial} faces of a simplex $\Delta$ come in two flavours: finitary faces and faces at
infinity.

More precisely, given $\ul{X}\in\Og(\Gamma)$, a {\em finitary face} of $\Delta_X$ corresponds to 
the collapse a forest in $\ul X$ whose components have trivial core, so that
the resulting graph of groups induces the same splitting of the fundamental group (as graph of
groups).
  We denote the finitary faces just {\em faces}. We define the 
   {\em closed} simplex $\overline{\Delta}$  as the closure of $\Delta$ in $\O(\Gamma)$,
   that is: $$\overline{\Delta}=\Delta\cup\{\text{all
     faces of $\Delta$}\}=\Delta\cup\{\text{all finitary faces of $\Delta$}\}.$$ 

The  {\em finitary boundary} of $X$ is the set of its
  proper finitary faces: $$\partial_\O\Delta=\partial_\O\overline\Delta=\overline\Delta\setminus\Delta.$$

A non-finitary simplicial face of an open simplex $\Delta_X$, corresponds to the
collapse of sub-graph $\ul A\subset \ul X$ with at least a component with non-trivial core, and    
belongs to the outer space $\O(X/A)$, (instead of $\O(X)$).  However, if $Y=X/A$, the simplicial
topology naturally defines a topology on $\overline{\Delta_X}\cup \overline{\Delta_Y}$, which
we still call the simplicial topology. Such a simplicial face will be called a {\em face at
  infinity of $\overline{\Delta_X}$}, and if all components of $\ul A$ are core-graphs, we call
it a {\em face at infinity of $\Delta_X$}. So, with this terminology, any simplicial face of
$\Delta$ is either a finitary face of $\Delta$, or a face at infinity of some finitary (not
necessarily proper) face of $\Delta$. We refer to~\cite[Section 2]{FMpartI} for a more detailed discussion.

We define the {\em boundaries at
  infinity} of a simplex $\Delta$ by 
$$\partial_\infty\Delta=\{\text{faces at infinity of } \Delta\}\qquad\text{(collapsing of only core sub-graphs)}$$
$$\partial_\infty\overline\Delta=\{\text{faces at infinity of } \overline\Delta\}
\quad\text{(more general collapsing)}$$

and the {\em closure at infinity by}
$$\overline\Delta^\infty=\overline\Delta\cup\partial_\infty\overline\Delta.$$

If we denote by $\partial \Delta$ the simplicial boundary of $\Delta$, we have
$$\partial \Delta=\partial_\infty\overline \Delta\cup\partial_\O\overline\Delta$$
and $$\partial_\infty \overline\Delta=\bigcup_{F=\text{face of }\Delta}\partial_\infty F$$
(where the union is over all faces of $\Delta$, $\Delta$ included.) Moreover, the simplicial
closure of $\Delta$ is just $\overline\Delta^\infty$.

  We define the {\em boundary at infinity} and the {\em simplicial bordification} of $\O(\Gamma)$ as
$$\partial_\infty\O(\Gamma)=\bigcup_{\Delta\text{ simplex}}\partial_\infty\Delta
\qquad \text{and}\qquad \overline{\O(\Gamma)}=\overline{\O(\Gamma)}^\infty=\O(\Gamma)\cup\partial_\infty\O(\Gamma).$$

\begin{rem}
	We note that when $\Gamma$ is just a topological graph with $\pi_1(\Gamma)=F_n$ (all
        vertex groups are trivial)  then $\O(\Gamma)$ is simply the Culler-Vogtmann Outer
        space $CV_n$, and the bordification $\overline{\O(\Gamma)}$ is the free splitting
        complex  $\mathcal{FS}_n$. (See~\cite{FMpartI} for more details).
\end{rem}

\subsection{Horoballs and regeneration (\cite[Section 2.7]{FMpartI})}\label{sechor}
We keep Notation~\ref{not:gamma}.
\begin{defn}[Horoballs]
Given $X\in\partial_\infty\O(\Gamma)$, $\Hor(X)$ is the set of marked metric trees
$Y\in\O(\Gamma)$ such that $\ul X$ is obtained from $\ul Y$ by collapsing a proper family of
core sub-graphs. By convention, we set $\Hor(X)=X$ when $X\in \O(\Gamma)$ (and we use $\Hor(\ul
X)$ for graphs). In other words, $Y\in\Hor(X)$ if $\ul X$ is obtained by setting to zero the
edge-lengths of a proper family of core sub-graphs (note that this implies that $\Delta_X$ is a simplicial face of $\Delta_Y$). 
\end{defn}

$\Hor(X)$ can be regenerated from $X$ as follows\footnote{Lemma~\ref{horconnected}, even if
  implicitly contained and proved in~\cite[Section 2.7]{FMpartI}, it is not explicitly stated
  there. We state and prove it here for future reference.}

\begin{lem}
	\label{horconnected}
	 Suppose $X\in\partial_\infty\Og(\Gamma)$. Let $Y\in\Og(\Gamma)$ and 
         $A=\sqcup_iA_i\subset Y$ be a family of core-graphs  such that $X=Y/A$.
         Then, for some $k_i$, we have
	
	$$\Hor(\wt{X})=\Pi_i\O(\wt{A_i},k_i).$$
	
	In particular, $\Hor(X)=\Hor(\wt{X})$ is path connected.
\end{lem}
\begin{rem}
	Note that we are using the tilde notation here, despite the objects being equivalent, to emphasise that the marked points are points in the trees.
\end{rem}
\proof 

Let
$v_i$ be the non-free vertex of $X$ corresponding to $A_i$. In order to recover a generic point
$Z\in\Hor(X)$, we need to replace each $v_i$ with an element $V_i\in\Og(A_i)$. Moreover,
in order to completely define the marking on $Z$, we need to know where to attach  - to $V_i$ -  the edges of $X$
incident to $v_i$, and this choice has to be done in the universal covers $\wt{V_i}$. No more
is needed. Therefore, if $k_i$ denotes the valence of the vertex $v_i$ in $X$, we have
$$\Hor(\wt{X})=\Pi_i\O(\wt{A_i},k_i).$$
(Note that some $k_i$ could be zero, e.g.  if $A_i$ is a connected component of $Y$.)

Each of the spaces $\O(\wt{A_i},k_i)$ is path connected. Indeed, the map that `forgets' the marked points is a continuous map to a path connected space whose fibers are connected; since each $A_i$ is connected, we can continuously deform any marked $k$-tuple of points to another, as we do not insist that they are distinct. 

The last statement now follows since a product of path connected spaces is path connected. 
\qed

\medskip

\begin{rem}
\label{defining pi}	

With above notation, the forgetting of marked points, gives a well-defined projection $\Hor(X)\to
\O(A)=\Pi_i\O(A_i)$. In what follows we will be mainly interested in the composition of such
map with the projection $\O(A)\to\p\O(A)$. We therefore give a name to such projection, defining

$$\pi:\Hor(X)\to\mathbb P\O(A).$$

(Here $\Hor(X)$ is intended to be not projectivized).
\end{rem}

\begin{rem}
Note that {\em the same} tree $X$ can be considered as a point at
infinity of {\em different} spaces. If we need
to specify in which space we work we write $\Hor_\Gamma(X)$.
\end{rem}

\subsection{Displacement function, optimal maps and train tracks}\label{s3.4}
For any $g\in G$ and $X\in \O(\G)$, the translation length $L_X(g)$ is defined as $\inf_p=d_X(gp,p)$. Elements with zero translation length correspond to vertex stabilisers, and are
called {\em elliptic}; others have the infimum realised along an axis, and are referred to as
{\em hyperbolic} elements. (Note that an element being elliptic or hyperbolic depends only on
$\G$ and not on $X$). The same happens in $\O(\Gamma)$ componentwise (that is for
$g\in\cup_i H_i$, where $H_i$ is as in Notation~\ref{not:gamma}). In this section we consider only hyperbolic elements. 

Given $X,Y\in O(\Gamma)$, we can compute the translation length of any hyperbolic $g\in \cup_i H_i$ in both
$X$ and $Y$, and we define 
$$\Lambda(X,Y)=\sup_g\frac{L_Y(g)}{L_X(g)}=
\inf\{\Lip(f) \ : \  f:X\to Y \text{ Lipschitz equivariant map}\}$$
where $\Lip(f)$ denotes the best Lipschitz constant for $f$.

It turns out that above second inequality is indeed true, and that $\sup$ and $\inf$ are in fact $\max$ and $\min$
(Theorem~\ref{sausagelemma}(\cite[Theorem~3.7]{FMpartI}), and
Theorem~\ref{Lemma_opt}(\cite[Theorem~3.15]{FMpartI})). $\Lambda(X,Y)$ can be computed by means of {\bf straight maps}; that is to say 
equivariant Lipschitz maps with constant speeds on edges.
Given a straight map, the {\bf tension graph} $X_{\max}(f)$ (or simply
$X_{\max}$) is the union of edges that are maximally stretched by $f$. A straight map that
realises the above minimum is called {\bf weakly optimal map}, and it is {\bf optimal} if the
tension graph has no one-gated vertex (we refer to~\cite{FMpartI} for further details on gate-structures). An optimal map is {\bf minimal}
if the tension graph coincides with the union of the axes of all maximally stretched elements. 

\begin{rem}
	\label{distance on boundary}
	We could also take the following point of view: given $X, Y \in \overline{\O(\Gamma)}$, let $Hyp(X)$ denote the set of hyperbolic elements in $X$, and similarly for $Hyp(Y)$. Note that if $X, Y \in \O(\Gamma)$, then $Hyp(X) = Hyp(Y)$. One can then define, 
	$$
	\Lambda(X,Y) = \sup_{g \in Hyp(X)} \frac{L_Y(g)}{L_X(g)} = \inf\{\Lip(f) \ : \  f:X\to Y \text{ Lipschitz equivariant map}\}$$
	where $\Lip(f)$ denotes the best Lipschitz constant for $f$, as long as $Hyp(Y) \subseteq Hyp(X)$. That is, as long as elliptic elements of $X$ are also elliptic in $Y$. If this is not the case, we set $\Lambda(X,Y) = \infty$. 
\end{rem}

  For any automorphism $[\phi] \in\Out(\Gamma)$ we define the displacement function
$$\lambda_\phi:\O(\Gamma)\to\R\qquad \lambda_\phi(X)=\Lambda(X,\phi X)$$
If $\Delta$ is a simplex of $\O(\Gamma)$ we define
$$\lambda_\phi(\Delta)=\inf_{X\in\Delta}\lambda_\phi(X)$$
If there is no ambiguity we write simply $\lambda$ instead of $\lambda_\phi$.
Finally, we set
$$\lambda(\phi)=\inf_{X\in\O(\Gamma)}\lambda_\phi(X)$$

In~\cite{FMpartI} the behaviour of the displacement near points in
$\partial_\infty(\O(\Gamma))$ is extensively studied. In particular, it is proven that if 
$X_\infty\in\partial_\infty(\O(\Gamma))$ is the limit of a sequence of points
$X_i\in\O(\Gamma)$ such that $\lambda_\phi(X_i)$ is bounded above, then $X_\infty$ and $\phi
X_\infty$ have the same elliptic elements, and $\phi$ induces an element of
$\Out(X_\infty)$. Therefore, for those points, the expression  $\lambda_\phi(X_\infty)=\Lambda(X_\infty,\phi X_\infty)$
still makes sense. For other points $T\in \partial_\infty(\O(\Gamma))$ we set $\lambda_\phi(T)=\infty$.

\begin{rem}
\label{displacement at infinity}\label{reminf}
Observe that $\Lambda(X_\infty, \phi(X_\infty))$ is finite, according to Remark~\ref{distance on boundary}, if and only if the set of elliptic elements of $X_\infty$ is $\phi$-invariant. If $X_\infty \in \partial_\infty(\O(\Gamma))$ has finite $\phi$-displacement, then we can regenerate $X_\infty$ to a point $X \in \O(\Gamma)$ such that $X$ admits an invariant core subgraph, $A$, which (as a forest) is a union of the axes of the elliptic elements of $X_\infty$ which are hyperbolic in $X$. $X_\infty$ is obtained from $X$ by collapsing $A$. Then, by Lemma~\ref{lemma9}, there will be a sequence $X_i\in\O(\Gamma)$ such that $X_i \to X$ and $\lambda_\phi(X_i)$ is bounded above.

Moreover, if $[\phi]$ is irreducible, then \textit{every } $X \in \partial_\infty(\O(\Gamma))$ has infinite $\phi$-displacement since no point in $\O(\Gamma)$ admits an invariant core graph. 
\end{rem}

The displacement function of an automorphism is not continuous at the bordification. Given
$[\phi]$, we say
that $X\in\overline{\O(\Gamma)}$ has {\bf not jumped} if there is a sequence $X_i\to X$ of points in
$\O(\Gamma)$ such that $\lambda_\phi(X_i)\to\lambda_\phi(X)$. Given a simplex $\Delta$ with $X$ in the
boundary at infinity of $\Delta$, we say that $X$ has {\bf not jumped in $\Delta$} if the above
condition holds with $X_i\in\Delta$.

\begin{defn}[$\O$-maps]\label{defOmap} For $X,Y\in O(\G)$, a map $f:X\to Y$ is called $\O$-map
  if it is Lipschitz continuous and $G$-equivariant. 
	Let now $X=(X_1,\dots,X_k)$ and $Y=(Y_1,\dots,Y_k)$ be two elements of $\O(\Gamma)$.
	A map $f=(f_1,\dots,f_k):X\to Y$ is called an $\O$-map if for each $i$ the map
	$f_i$ is an $\O$-map from $X_i$ to $Y_i$. (No index permutation here). 
\end{defn}

Let $[\phi]=[(\sigma,\phi_1,\dots,\phi_k)]$ be an element of $\Out(\Gamma)$ - see Definition~\ref{autgamma}. 

\begin{defn}[Maps representing {$[\phi]$}]\label{maprepf}
  Let $X\in\O(\Gamma)$. We say that a map $f:X\to X$ {\em represents
    $[\phi]$}\footnote{In~\cite{FMpartI} we used {\em $f$ represents $\phi$}. Such notation
    appears in Section~\ref{s4} where we quote results from~\cite{FMpartI}.} if
  $f$ maps each $X_i$ to $X_{\sigma_{(i)}}$, and such that,  if we denote by $f_i$ the map
  $f|_{X_i}:X_i\to  X_{\sigma_{(i)}}$, then $f_i$ is Lipschitz and equivariant with respect to the isomorphism
  $\phi_i:H_i\to H_{\sigma_{(i)}}$, that is $f_i(hx)=\phi_i(h)f_i(x)$. Note that a map
  representing $[\phi]$ can be viewed as an $\O$-map $f:X\to \phi X$.
    We say that $f$ is optimal if each $f_i$ is optimal.
	
	If $X$ is a $\Gamma$-graph, then a map $f:X\to X$ represents $[\phi]$ if it has a lift 
	$\wt f:\wt X\to \wt X$ representing $[\phi]$.
\end{defn}

\begin{nt}
For notational coherence with~\cite{FMpartI},  if not otherwise specified, if
$X,Y\in\O(\Gamma)$ and $f:X\to Y$, when we say that $f$ is straight we
  understand that it is also an $\O$-map. 
\end{nt}

In~\cite[Section 4]{FMpartI} we introduced the notion of partial train tracks and partial train tracks at
infinity. Roughly: given $[\phi]$, a {\bf partial train track} for $[\phi]$ on $X\in
\O(\Gamma)$ is a straight map
$f:X\to X$ representing $[\phi]$ such that $X$ has a $f$-invariant sub-graph to which the
restriction of $f$ is a train track; a {\bf partial train track at infinity} is  when
$X\in\partial_\infty(\O(\Gamma))$. 

The deep link between partial train track maps and displacement function, is fully studied and exploited
in~\cite{FMpartI}. In this paper we use results from~\cite{FMpartI}, but we don't need to
directly involve partial train tracks. And in fact the words ``train track'' will appear only
in Section~\ref{s4}, where we quote literally statements form~\cite{FMpartI}.

For completeness of exposition we just recall that, as proved in~\cite{FMpartI},
given $[\phi]$, the minimally displaced set of $[\phi]$, that is to say the set of trees $T$
such that $\lambda_\phi(T)=\lambda(\phi)$,  coincides
with the set of points admitting a partial train track map. For reducible automorphisms, the minimally displaced set may be
empty in $\O(\Gamma)$, but if one includes partial train tracks
at infinity (partial train tracks for a point at the bordification where the displacement does
not jump)  then the set of points admitting these partial train tracks is always non-empty and is
contained in the minimally displaced set (of points at infinity).  Notationally, $\Min(\phi)=\TT(\phi)$
is the minimally displaced set in $\O(\Gamma)$, which coincides with the set of points supporting a partial train track map.

\section{Results needed from~\cite{FMpartI}}\label{s4}
In what follows, we will need to quote many lemmas and results
from~\cite{FMpartI}. For the ease of the reader we collected the statements we need
from~\cite{FMpartI} in this section.  We decided to quote them exactly as they appear
in~\cite{FMpartI}, paying the price that some of them may look a little
redundant or overstated here.  The reader can safely skip this section now, coming back
here when a needed result is cited.

\begin{thm}[Sausage Lemma~{\cite[Theorem 3.7]{FMpartI}}]\label{sausagelemma}
  Let $X,Y\in\Og(\Gamma)$. The stretching factor $\Lambda(X,Y)$ is realized by a loop
  $\gamma\subset X$ having  has one of the following forms:
  \begin{itemize}
  \item Embedded simple loop $O$;
  \item embedded ``infinity''-loop $\infty$;
  \item embedded barbel $O$--- $O$;
  \item singly degenerate barbel $\bullet$---$O$;
  \item doubly degenerate barbel $\bullet$---$\bullet$.
  \end{itemize}
(the $\bullet$ stands for a non-free vertex.) Such loops are usually named ``candidates''.
\end{thm}

\begin{thm}[{\cite[Theorem 3.15]{FMpartI}}]\label{Lemma_opt}
Let $X,Y\in\O(\Gamma)$ and let $f:X\to Y$ be a straight map. There is a map\footnote{We describe
  an algorithm to find the map $\wopt(f)$, but the algorithm will depend on some choices, hence
  the map $\wopt(f)$ may be not unique in general.} $\wopt(f):X\to Y$
which is weakly optimal and such that
$$d_\infty(f,\wopt(f))\leq \vol(X)(\lambda(f)-\Lambda(X,Y))$$
Moreover, for any weakly optimal map $\f:X\to Y$ and for any 
$\varepsilon >0$ there is an optimal map $g:X\to Y$ such that
$d_\infty(g,\f)<\varepsilon$.
\end{thm}

\begin{defn}[Exit points, {\cite[Definition 4.19]{FMpartI}}]\label{exitp}
    Let $[\phi] \in\Out(\Gamma)$. A point $X\in\O(\Gamma)$ is called an {\em exit point} of
    $\Delta_X$ if for any neighbourhood $U$ of $X$ in $\O(\Gamma)$ there is an optimal map
    $f:X\to X$, representing $\phi$, a point $X_E\in U$, and a folding path (\cite[Definition
    3.21]{FMpartI}) directed by $f$,  
    $X=X_0,X_1,\dots,X_m=X_E$ in $U$, such that $\Delta_{X_i}$ is a
    finitary face of $\Delta_{X_{i+1}}$, $\Delta_{X}$ is a proper face of
    $\Delta_{X_E}$, and such that $$\lambda_\phi(X_E)<\lambda_\phi(X) \qquad \text{(a strict inequality)}.$$
\end{defn}

\begin{lem}[{\cite[Lemma 4.20]{FMpartI}}]\label{LemmaX}
  Let $[\phi] \in\Out(\Gamma)$ and $X\in\O(\Gamma)$ such that
  $\lambda_\phi(X)$ is a local minimum for $\lambda_\phi$ in $\Delta_X$.
  Suppose $X\notin\TT(\phi)$.

Then, for any open neighbourhood $U$ of $X$ in $\O(\Gamma)$,
there is an optimal map $f:X\to X$, representing $\phi$, points  $Z,X'\in U$, and a folding
path, $X=X_0,\dots,X_m=Z,X_{m+1},\dots,X_n=X'$, directed by $f$ and such that: 
\begin{itemize}
\item $X_0,\dots,X_m\in U\cap \Delta_X$,
\item $\lambda_\phi(Z)=\lambda_\phi(X)$,
\item $\Delta_X$ is a proper face of $\Delta_{X'}$,
\item $\lambda_\phi(X')<\lambda_\phi(X)$.
\end{itemize}
In particular $X$ is an exit point of $\Delta_X$.
\end{lem}

\begin{thm}[{\cite[Theorem 5.8]{FMpartI}}, lower semicontinuity of $\lambda$]\label{fatto1}
  Fix $\phi\in\Aut(\Gamma)$ and $X\in\Og(\Gamma)$. Let $(X_i)_{i\in\N}\subset \Delta_X$ be a
  sequence such that there is $C$ such that for any $i$,  $\lambda_\phi(X_i)<C$. Suppose that $X_i\to
  X_\infty\in\partial_\infty\Delta_X$ which is  obtained from $X$  by collapsing a sub-graph
  $A\subset X$. Then $\phi$ induces an element of $\Aut(X/A)$, still denoted by
  $\phi$.

  Moreover $\lambda_\phi(X_\infty)\leq \liminf_{i\to\infty} \lambda_\phi(X_i)$, and
  if strict inequality holds, then there is a sequence of
  minimal optimal maps $f_i:X_i\to X_i$ representing $\phi$ such that eventually on $i$ we have
  $(X_i)_{\max}\subseteq \core(A)$\footnote{By~\cite[Proposition 5.6]{FMpartI} we know that
    $\core(A)$ is $\phi$-invariant}.
\end{thm}

\begin{lem}[{\cite[Lemma 5.12]{FMpartI}}, regeneration of optimal maps]\label{lemma9}
  Fix $\phi\in\Aut(\Gamma)$ and $X\in\Og(\Gamma)$. Let $X_\infty\in\partial_\infty\Delta_X$ be
  obtained from $X$ by collapsing a $\phi$-invariant core sub-graph $A$.
  Then, for any straight map $f_A:A\to A$ representing $\phi|_A$, and for any $\e>0$
  there is $X_\e\in\Delta_X$ such that $$\lambda_\phi(X_\e)\leq
  \max\{\lambda_\phi(X_\infty)+\e, \Lip(f_A)\}.$$

More precisely, for any $Y\in \mathbb P\Og(A)$ and map $f_Y:Y\to Y$ representing $\phi|_A$, for
any map $f:X_\infty\to X_\infty$ representing $\phi$,
for any $\widehat X\in\pi^{-1}(Y)$ \footnote{See Remark~\ref{defining pi} for an explanation of the map $\pi$.}, and for any $\e>0$; there is
$0<\delta=\delta(f,f_Y,X_\infty,\Delta_{\widehat X})$,
such that for any  $Z\in\Delta_{\widehat X}\cap\pi^{-1}(Y)$, if $\vol_Z(Y)<\delta$  there is a
straight map $f_Z:Z\to Z$
representing $\phi$ such that $f_Z=f_Y$ on $Y$ and
$$\Lip(f_Z)\leq\max\{\lambda_\phi(X_\infty)+\e, \Lip(f_Y)\}$$ (hence the optimal map
  $\opt(f_Z)$ satisfies the same inequality\footnote{We notice that while $f_Z=f_Y$ on $Y$,
  this is no longer true for $\opt (f_Z)$}).
\end{lem}

\begin{thm}[{\cite[Corollary 5.14]{FMpartI}}]\label{newjump}
  Let $\phi\in\Aut(\Gamma)$. Let $X\in\Og(\Gamma)$ containing an invariant
  sub-graph $A$. Let $X_\infty=X/A$ and $C=\core(A)$. Then
  $$\lambda_{\phi|_C}(\Delta_C)\leq\lambda_\phi(\Delta_X).$$ Moreover the following are
  equivalent:
  \begin{enumerate}
  \item $X_\infty$ has not jumped in $\Delta_X$;
  \item $\lambda_\phi(X_\infty)\geq\lambda_{\phi}(\Delta_X)$;
  \item $\lambda_\phi(X_\infty)\geq\lambda_{\phi|_C}(\Delta_C)$.
  \end{enumerate}
 In particular, $\lambda_\phi(X_\infty)$ cannot belong to the (potentially empty)  interval
$(\lambda_{\phi|_C}(\Delta_C),\lambda_\phi(\Delta_X))$. Moreover, points realising
$\lambda_\phi(\Delta_X)$ do not jump in $\Delta_X$.
\end{thm}

\begin{cor}[{\cite[Corollary 5.17]{FMpartI}}]\label{fatto2}
  Let $\phi\in\Aut(\Gamma)$. Let $\Delta$ be a simplex of $\Og(\Gamma)$. 
  Then there is a min-point $X_{\min}$ in $\overline{\Delta}^\infty$ (\em{i.e.} a point so that
  $\lambda_\phi(X_{\min})=\lambda_\phi(\Delta)$; note that $X_{\min}$ does not jump in $\Delta$ by
  Theorem~\ref{newjump}). 

Moreover, suppose that $X_{\min}$ is {\em maximal} in the following sense: if $X' \in \overline{\Delta}^\infty$ such that $\lambda_\phi(X') =\lambda_\phi(X_{\min})=\lambda_\phi(\Delta)$, and $\Delta_{X_{\min}} \subseteq \overline{\Delta_{X'}}^\infty$, then $\Delta_{X_{\min}} = \Delta_{X'}$. ($X_{\min}$ is maximal with respect to the partial order induced by the faces of $\Delta$). Then: 
	\begin{itemize}
        \item $\lambda_{\phi}(X_{\min}) = \lambda_{\phi}(\Delta_{X_{\min}}) =
        \lambda_{\phi}(\Delta)$;
      \item any point $P$, such that $\Delta_{X_{\min }} \subseteq \overline{\Delta_P}^\infty \subseteq
        \overline{\Delta}^\infty$, satisfies $\lambda_{\phi}(P) \geq \lambda_\phi(\Delta)$ (hence does not
        jump in $\Delta$ by Theorem~\ref{newjump});  
        \item for any $\epsilon >0$, there exist points $Z, W$ such that:
          \begin{itemize}
		\item $Z \in \Delta$,
		\item $\Delta_{X_{\min}} \subseteq \overline{\Delta_W}^\infty \subseteq
                  \overline{\Delta}^\infty$,  
		\item $\lambda_{\phi}(W), \lambda_{\phi}(Z) \leq \lambda_{\phi}(\Delta) + \epsilon$,
		\item $\lambda_{\phi}$ is continuous along the Euclidean segments, $ZW$ and
                  $WX_{\min}$, and any point $P$ along these segments satisfies the following:
                  $\lambda_{\phi}(\Delta) \leq \lambda_{\phi}(P)$.
          \end{itemize}
		
	\end{itemize}

(We allow degeneracies, meaning that $X_{min}$ could equal $W$, or even $Z$).
\end{cor}

\begin{lem}[{\cite[Lemma 6.2]{FMpartI}}]\label{lconvexity}
  For any $[\phi]\in\Out(\Gamma)$ and for any open simplex $\Delta$ in $\O(\Gamma)$ the function
  $\lambda=\lambda_\phi$ is quasi-convex\footnote{A function $f:[a,b]\to\R$
    is quasi-convex if for any $a\leq x\leq t\leq y\leq b$ we have
    $f(t)\leq\sup\{f(x),f(y)\}$.} on segments of $\Delta$. Moreover, for $A,B \in \Delta$, if
  $\lambda(A)>\lambda(B)$ then  $\lambda$ is strictly monotone near $A$\footnote{In this
    statement $A,B$ are points of $\Delta$, and monotonicity is referred to the restriction of
    $\lambda$ to the segment joining $A,B$.}.
\end{lem}

\begin{lem}[{\cite[Lemma 6.3]{FMpartI}}]\label{lconv2}
  Let $[\phi]\in\Out(\Gamma)$, let $\lambda=\lambda_\phi$, and let $\Delta$ be a simplex in
  $\O(\Gamma)$. Let $A,B\in \overline{\Delta}^\infty$ be two points that have not jumped in
  $\Delta$. Then for any 
  $P\in\overline{AB}$ $$\lambda(P)\leq\max\{\lambda(A),\lambda(B)\}$$

  Moreover, if $\lambda(A)\geq \lambda(B)$, then $\lambda|_{\overline{AB}}$ is continuous at $A$.
\end{lem}
\begin{thm}[{\cite[Theorem 7.2]{FMpartI}}]\label{conj}
For any $\Gamma$ the global simplex-displacement spectrum
$$\operatorname{spec}(\Gamma)= \Big\{\lambda_\phi(\Delta): [\phi]\in\Out(\Gamma),
\Delta\text{ a simplex of } \overline{\O(\Gamma)}^\infty \text{such that }\lambda_\phi(\Delta)<+\infty \}$$ is well-ordered as a subset of $\mathbb R$.
In particular, for any $[\phi]\in\Out(\Gamma)$ the spectrum of possible minimal displacements $$\operatorname{spec}(\phi)= \Big\{\lambda_\phi(\Delta):\Delta\text{ a simplex of }
\overline{\O(\Gamma)}^\infty \text{ such that } \lambda_\phi(\Delta)<+\infty\}$$ is well-ordered as a subset of $\mathbb R$.
\end{thm}

\begin{thm}[{\cite[Theorem 7.3]{FMpartI}}]\label{thmminptE} Let $\Gamma$ be as in Notation~\ref{not:gamma}.
  Let $[\phi]$ be any element in $\Out(\Gamma)$. Then there exists $X\in\overline{\O(\Gamma)}^\infty$
  that has not jumped and such that $$\lambda_\phi(X)=\lambda(\phi).$$
\end{thm}

\begin{lem}[{\cite[Lemma 7.7]{FMpartI}}]\label{dom5}
  Let $\phi\in\Aut(\Gamma)$. Let $X_\infty\in\overline{\Og(\Gamma)}$ which has not jumped.
  Suppose that there is a loop $\gamma \in X_\infty$ and $k>0$ and such that
  $L_{X_\infty}(\phi^n (\gamma)) \geq k^n L_{X_\infty}(\gamma)$ for all $n \in \N$. 
Then,  $$k\leq\lambda(\phi).$$

  In particular, if $X_\infty$ is a train track for $\phi$ as an element of $\Aut(X_\infty)$,
  then it is a minpoint  for $\phi$ as an element of $\Aut(\Gamma)$.
\end{lem}

\begin{thm}[{\cite[Theorem 7.8]{FMpartI}}]\label{corlalx}
  Let $\phi\in\Aut(\Gamma)$. Let $X\in\O(\Gamma)$ and $X_\infty$ be such that $\ul{X_\infty}$ is obtained from $\ul{X}$ by
  collapsing a $\phi$-invariant core sub-graph $\ul{A}$. Then $$\lambda(\phi|_A)\leq \lambda(\phi).$$

Moreover, if $\lambda(\phi|_A)=\lambda_\phi(X_\infty)$, then $$\lambda(\phi)=\lambda(\phi|_A).$$

In particular $X_\infty$ has not jumped
if and only if $$\lambda(\phi)\leq \lambda(X_\infty).$$
\end{thm}

\begin{rem}
	We note that if a point has not jumped, this simply means that there is some sequence converging to it, whose displacements tend to the displacement of that point. In general this will not hold for all sequences tending to the point.  
\end{rem}

\begin{thm}[{\cite[Theorem 7.13]{FMpartI}}]\label{strongcorred}
  Let $\phi\in\Aut(\Gamma)$. Let $X$ be a $\Gamma$-graph having a $\phi$-invariant core
  sub-graph $A$. Then there is $Z\in\overline{\O(X/A)}^\infty$ and $W\in\Hor_{\O(\Gamma)}(Z)$
  such that the simplex $\Delta_W$ 
  contains a minimising sequence for $\lambda$. Moreover if $Y$ is the graph used to
  regenerate $W$ from $Z$, then the minimising sequence can be chosen with straight maps $f_i$
  such that $f_i(Y)=Y$ and $\Lip(f_i)\to\lambda(\phi)$.
\end{thm}



\section{Statement of the connectedness theorem and regeneration of paths in the bordification}
We recall here Notation~\ref{not:gamma} (as a courtesy for readers who skipped the
first sections).

\begin{itemize}
\item  $\Gamma=\sqcup \Gamma_i$ will always mean that $\Gamma$ is a finite disjoint
union of finite graphs of groups $\Gamma_i$, each with trivial edge-groups and non-trivial fundamental group $H_i=\pi_i(\Gamma_i)$, each
$H_i$ being equipped with the splitting given by the vertex-groups. We set
$\rank(\Gamma)=\sum_i\rank(\Gamma_i).$
\end{itemize}

 We also recall that for any $[\phi] \in\Out(\Gamma)$ we defined the displacement function
$$\lambda_\phi:\O(\Gamma)\to\R\qquad \lambda_\phi(X)=\Lambda(X,\phi X)$$
If $\Delta$ is a simplex of $\O(\Gamma)$ we define
$$\lambda_\phi(\Delta)=\inf_{X\in\Delta}\lambda_\phi(X).$$
If there is no ambiguity we write simply $\lambda$ instead of $\lambda_\phi$.
Finally, we set
$$\lambda(\phi)=\inf_{X\in\O(\Gamma)}\lambda_\phi(X)$$
By convention (see Section~\ref{s3.4}) we extend the function $\lambda$ to points in
$X_\infty\in\partial_\infty(\O(\Gamma))$ for which there is a sequence of points
$X_i\in\O(\Gamma)$ such that $X_i\to X_\infty$ with $\lambda(X_i)$ bounded above, and we set
$\lambda=\infty$ on other points.

Finally, we recall that outer space comes in two flavours: trees and
graphs. We will chose which one we use on a case-by-case basis, depending on which is more convenient. For that purpose we introduced the notation
``$\O(\Gamma)$'' for trees and ``$\Og(\Gamma)$'' for graphs. Clearly $\Og(\Gamma)$ and
$\O(\Gamma)$ are isomorphic via $X\leftrightarrow \wt X$, and thus in all statements they are
completely interchangeable.

\begin{defn}\label{def:sp}
  Let $X,Y\in\overline{\O(\Gamma)}^\infty$. A {\em simplicial path} $\Sigma$ between $X,Y$ is
  given by: 
  \begin{enumerate}
  \item A finite sequence of points $X=X_0,X_1,\dots,X_k=Y$, called vertices, such that
    $\forall i=1,\dots, k$,  there is a simplex $\Delta_i$ such that
    $\Delta_{X_{i-1}}$ and $\Delta_{X_i}$ are  both simplicial\footnote{We
    remind that simplicial faces include faces at infinity. That is to say,
    $\Delta_{X_{i-1}}$ and $\Delta_{X_i}$ are both faces of $\overline{\Delta_i}^{\infty}$.}
  faces of $\Delta_i$. We allow 
    one of $\Delta_{X_{i-1}},\Delta_{X_i}$, or even both, to coincide with $\Delta_i$.
  \item Euclidean segments $\overline{X_{i-1}X_i}\subseteq \overline{\Delta_i}^\infty$, called
    edges. We require the interior of $\overline{X_{i-1}X_i}$
    (i.e. $\overline{X_{i-1}X_i}\setminus\{X_{i-1},X_i\}$) to be contained in $\Delta_i$.
   \item We say that $\Sigma$ is {\em alternating} if for every $i$ either $\Delta_{X_i}$ is a
     simplicial face of $\Delta_{X_{i-1}}$ or $\Delta_{X_{i-1}}$ is a simplicial face of
     $\Delta_{X_i}$. Note that any simplicial path can be made alternating just by adding some
     extra vertex. 
 \end{enumerate}
\end{defn}

\begin{defn}
  We say that a set $\chi$ is {\em connected by simplicial paths} if for any $x,y\in\chi$ there is a
  simplicial path between $x$ and $y$ which is entirely  contained in $\chi$. 
\end{defn}

\begin{thm}[Level sets are connected]\label{tconnected}
Let $[\phi] \in\Out(\Gamma)$. For any $\e > 0 $ the set
$$\{X\in\O(\Gamma):\lambda_\phi(X)\leq \lambda(\phi) + \e \}$$ is connected in $\O(\Gamma)$ by simplicial paths.
The set
$$\{X\in\overline{\O(\Gamma)}^\infty:\lambda_\phi(X)=\lambda(\phi)\}$$
is connected by simplicial paths in $\overline{\O(\Gamma)}^\infty$.

Moreover, connecting paths can be chosen so that the displacement $\lambda_\phi$ is continuous
along them.
\end{thm}
The main goal of the paper is the proof of Theorem~\ref{tconnected}. The rough
strategy is to prove that paths in the bordification can
regenerate to paths in $\O(\Gamma)$ without increasing $\lambda$ too much. Then, 
the first claim will follow from the second, which we will prove via a peak-reduction
argument. Proofs proceed via induction on the rank of $\Gamma$. This is part of the reason
that we need to fundamentally deal with the case where $\Gamma$ is disconnected. We remind
that Theorem~\ref{tconnected}, if specialised to the case where $\Gamma$ is connected and
vertex-groups are trivial, is a $CV_n$-statement about connectedness of level sets of, not
necessarily irreducible, automorphisms of the free group $F_n$.

\begin{cor}
	\label{tconnectedirred} 
	Let $[\phi] \in \Out(\Gamma)$ be irreducible. Then the set
	$$\{X\in\O(\Gamma):\lambda_\phi(X) = \lambda(\phi)  \}$$ is connected in $\O(\Gamma)$ by simplicial paths.
\end{cor}
\begin{proof}
	This is an immediate consequence of Theorem~\ref{tconnected}, since by Remark~\ref{displacement at infinity}, the irreducibility of $[\phi]$ implies that every point on the boundary, $\partial_\infty(\O(\Gamma))$, has infinite displacement.
\end{proof}

\begin{rem}
  Theorem~\ref{tconnected} is trivially true if $\rank(\Gamma)=1$, because in that case
  either $\O(\Gamma)$ or $\mathbb P\O(\Gamma)$ is a single point.
\end{rem}

\begin{lem}[Regeneration of segments]\label{regseg}
Fix $[\phi]\in\Out(\Gamma)$. Let $X_\infty,Y_\infty\in\overline{\O(\Gamma)}^\infty$ such that
$\Delta_{Y_\infty}$ is a (not necessarily proper) simplicial face of
$\Delta_{X_\infty}$. Suppose that  $\lambda(X_\infty)\geq\lambda(\phi)$.
Then there is an open simplex $\Delta$ of $\O(\Gamma)$ such
that for any $\e>0$ there is  $Y\in \Hor(Y_\infty)\cap \overline\Delta$ and
$X\in\Hor(X_\infty)\cap\Delta$
 such that\footnote{Note that the fact that $\Hor(Y_\infty)\cap\overline\Delta\neq\emptyset$
   implies that $\Delta_{Y_\infty}$ is a simplicial face of $\Delta$. The same holds true for $X_\infty$.} $$\lambda_\phi(Y),\lambda_\phi(X)<\max\{\lambda_\phi(Y_\infty),\lambda_\phi(X_\infty)\}+\e.$$
 Moreover, such an inequality holds for any  $T\in\overline{XX_\infty}$ and any $S\in \overline{YY_\infty}$.
\end{lem}
\proof For this proof will be more convenient to work in $\Og$ rather than $\O$. Let $X_\infty$ be obtained by collapsing a $\phi$-invariant core-subgraph $A$ from a
$\Gamma$-graph $\widehat X$. Since $\lambda_\phi(X_\infty)\geq\lambda(\phi)$, by Theorem~\ref{corlalx} $\lambda(\phi|_A)\leq\lambda_\phi(X_\infty)$.
By Theorems~\ref{thmminptE} and~\ref{strongcorred}, there is a simplex in $\Og(A)$ that contains a minimising sequence
for $\lambda(\phi|_A)$. Let $A_\e$ be a point in that simplex such that
$\lambda(A_\e)<\lambda(\phi|_A)+\e$. The required simplex $\Delta$ is
obtained by inserting a copy of $A_\e$ in place of $A$ in $X_\infty$. We note that such
a $\Delta$ is not unique.
 By Lemma~\ref{lemma9} there is a
point $X\in\Delta\cap\Hor(X_\infty)$ such that
$\lambda_\phi(X)\leq\lambda_\phi(X_\infty)+\e$.

Consider now the points in $\overline\Delta\cap \Hor(Y_\infty)$.
By hypothesis there is a $\phi$-invariant $B\subseteq X_\infty$ such that
as a graph (that is, forgetting the metric), $Y_\infty$ is obtained from $X_\infty$ by collapsing
$B$. $B$  has a pre-image in $X$ still denoted by $B$. Let $T$ be the forest $(A\cup
B)\setminus \core(A\cup B)$. If $Y'=X/T$, as a graph, $Y_\infty=X/(A\cup B)=Y'/\core(A\cup B)$.

Thus the finitary face $\Delta_{Y'}$ of $\Delta$ obtained by the collapse of $T$ intersects
$\Hor(Y_\infty)$.

Let $f:X\to X$ be an optimal map representing $[\phi]$. Since $\core(A\cup B)$ is $\phi$-invariant,
$f(\core (A\cup B))\subset \core(A\cup B)$ up to homotopy. It follows that there is
a straight map $g:\core(A\cup B)\to \core(A\cup B)$ representing $[\phi|_{A\cup B}]$ such that
$\Lip(g)\leq\lambda_\phi(X)\leq\lambda_\phi(X_\infty)+\e$.
By Lemma~\ref{lemma9} there is a point $Y\in\Hor(Y_\infty)\cap \Delta_{Y'}$ such that
$\lambda_\phi(Y)\leq\max\{\lambda_\phi(Y_\infty)+\e,\Lip(g)\}\leq\max\{\lambda_\phi(Y_\infty)+\e,\lambda_\phi(X_\infty)+\e\}$. The
last claim also follows by Lemma~\ref{lemma9}, since the volume of $A$ (or $B$) is strictly decreasing on the Euclidean segment $\overline{XX_\infty}$ (or $\overline{YY_\infty}$), and the invariant subgraph is being scaled uniformly. \qed

\medskip

Now we can plug in the inductive hypothesis in the proof of Theorem~\ref{tconnected}.

\begin{lem}[Regeneration of horoballs]\label{inhor}
   Suppose that Theorem~\ref{tconnected} is true in any rank less than $\rank(\Gamma)$. Let
  $[\phi] \in\Out(\Gamma)$. Let $T\in\Og(\Gamma)$ be a $\Gamma$-graph having a proper
  $\phi$-invariant core sub-graph
  $S$. Let $X\in\partial_\infty\Og(\Gamma)$ be the graph obtained from $T$ by collapsing $S$, and
  let $A,B\in\Hor(X)\subset\Og(\Gamma)$. Let $m_A$
  and $m_B$ be the supremum of $\lambda_\phi$ on the Euclidean segments $\overline{AX}$ and
  $\overline{BX}$ respectively. Then, for any $\e>0$ there is a simplicial path $\Sigma$
  between $A$ and $B$, and  in $\Hor(X)$, such that for any vertex $Z$ of $\Sigma$ we have $$\lambda_\phi(Z)<\max\{m_A,m_B\}+\e.$$
\end{lem}
\proof Let $L=\max\{m_A,m_B\}$. The displacement $\lambda_\phi(T)$ is a finite number just
because $\lambda_\phi$ is a well-defined function on $\Og(\Gamma)$. 
For any group element $g\in\cup_i H_i$, and for any $t\in[0,1)$, 
the translation length of $g$ in $T_t=tX+(1-t)T$ is bounded by $L_T(g)$.
Moreover, as $T$ is a finite graph, there is $\delta>0$ such that for any reduced loop $\gamma$ in
$T$, either $\gamma\subseteq S$ or the length of $\gamma$ in any $T_t$ is at least $\delta$.
Therefore, since $S$ is $\phi$-invariant, and by using Theorem~\ref{sausagelemma}, we see that 
$\lambda_\phi(T_t)$ is bounded by some constant $C$, uniform on $t$.

By Theorem~\ref{fatto1} we have
that $\lambda_\phi(X)$ is finite, and by Lemma~\ref{lemma9} both $m_A$ and $m_B$ are
finite\footnote{One has to apply Lemma~\ref{lemma9} as follows: $X$ here plays the role of
  $X_\infty$ of lemma; $T$ here is $X$ in lemma, $S$ here is $A$ in lemma; $A,B$ here play the
  role of $\widehat X$ in lemma (for suitable $Y$).}.
 
Recall that if $X=T/S$ as graphs of groups, then we  denote by $\pi:\Hor(X)\to \mathbb P\Og(S)$ the
projection that associates to a point in $\Hor(X)$ its collapsed part (see Section~\ref{sechor}).

For any $Y\in\Hor(X)$, $\lambda_\phi(Y)$ can be computed by tacking the supremum of stretching
factors of candidates given by Theorem~\ref{sausagelemma}. Those may or may not be contained in
$S$, and clearly the supremum over all candidates is bigger or equal to the supremum over
candidates contained in $S$.  Since $S$ is $\phi$-invariant, this implies that
$\lambda_\phi(\pi(Y))\leq\lambda_\phi(Y)$ for any $Y\in\Hor(X)$; so
$$\lambda_\phi(\pi(A))\leq\lambda_\phi(A)\quad
\lambda_\phi(\pi(B))\leq\lambda_\phi(B)
$$
hence, $\lambda_\phi(\pi(A)),\lambda_\phi(\pi(B))\leq L$. The rank of $S$ is strictly
smaller than $\rank(\Gamma)$ because it is a proper sub-graph of $T$. Hence Theorem~\ref{tconnected} holds for $\Og(S)$. Therefore, the induction hypothesis produces a finite simplicial path $(Y_i)\in\Og(S)$ between $\pi(A)$ and $\pi(B)$ such that $\lambda_\phi(Y_i)<L+\e$. Hence, by Lemma~\ref{horconnected}, there is
a finite  simplicial path in $\Hor(X)$ between $A$ and $B$
whose vertices are points $\widehat T_j$ such that for
any $j$ there is $i$ such that $\pi(\widehat T_j)=Y_i$. By Lemma~\ref{lemma9} there is a
simplicial path in $\Hor(X)$ whose vertices are points $Z_j\in\Delta_{\widehat T_j}$ such that $\pi(Z_j)=\pi(\widehat
T_j)=Y_i$ and $\lambda_\phi(Z_j)<L+\e$.\qed

\medskip

We recall that we are using the notation of Definition~\ref{def:sp}.
\begin{thm}[Regeneration of paths]\label{tregge}
  Suppose that Theorem~\ref{tconnected} is true in any rank less than $\rank(\Gamma)$. Let
  $[\phi] \in\Out(\Gamma)$. Let
  $\Sigma=(X_i)_{i=1}^m$ be an alternating simplicial path in $\overline{\O(\Gamma)}^\infty$,
  and let $L$ be a positive real number.

  Suppose that for any point $X_i$ we have
  $$\lambda(\phi)\leq\lambda_\phi(X_i)\leq L.$$

  Then, for any $\e>0$ there exists a simplicial path $\Sigma_\e = (W_j)_{j=1}^k$ in $\O(\Gamma)$, 
such that for any point $P$ of $\Sigma_\e$, $\lambda(P) \leq (L+\e)$.

Moreover, we can choose the path so that  $W_1\in\Hor(X_1)$, $W_k\in\Hor(X_m)$, each  $W_j$
belongs to the horoball of some $X_i$; and so that $X_1$ and  $X_m$  do not jump in
$\Delta_{W_1}$ and $\Delta_{W_k}$ respectively.
\end{thm}
\proof By Lemmas~\ref{lconvexity} and~\ref{lconv2}, and since the displacement is continuous in $\O(\Gamma)$, it
suffices to check displacement on vertices of $\Sigma_\e$.

For any $i<m$, we apply Lemma~\ref{regseg} to the $i^{th}$ pair of consecutive points
$X_i,X_{i+1}$. This produces points $A_i\in\Hor(X_i)$ and $B_{i+1}\in\Hor(X_{i+1})$ whose
displacement is less than $L+\e$. Note that $\e$ is arbitrary. In particular
Theorem~\ref{newjump} implies that $X_1$ does not jump in $\Delta_{A_1}$ and $X_m$ does not
jump in $\Delta_{B_m}$. Moreover, $A_i,B_{i+1}$ are in the same closed simplex of $\O(\Gamma)$
(so there is a Euclidean segment joining them). 

Additionally Lemma~\ref{regseg} tells us that the displacement of points along the segments, $\overline{A_i X_i}, \overline{B_i X_i}$ is bounded by $L+\e$.

Note that $A_i$ is defined for $1\leq i<m$ and $B_i$ for $1<i\leq m$. By Lemma~\ref{inhor}, for
$1<i<m$, there is a simplicial path $Y_{ij}$ between $B_i$ and $A_i$ such that
$Y_{ij}\in\Hor(X_i)$ and 
$\lambda_{\phi}(Y_{ij})\leq L+\e$. The path $\Sigma_\e$ is now defined by the concatenation of such
paths and the segments $\overline{A_iB_{i+1}}$. \qed

\section{Calibration of paths}
We keep Notation~\ref{not:gamma}. For the remaining of the section we fix
$[\phi]\in\Out(\Gamma)$.
We recall that for simplices $\Delta\subset\overline{\O(\Gamma)}^\infty$ we are using the notation $\lambda(\Delta)=\lambda_\phi(\Delta)=\inf_{X\in\Delta}\lambda_\phi(X)$.

Our aim is to run a peak reduction argument to prove Theorem~\ref{tconnected}, by starting with
a simplicial path and locally modifying it near peaks. Theorem~\ref{tregge} provides simplicial
paths with bounded displacement, however, for our purposes we need paths, that possibly touch
the boundary at infinity,  where the displacement
is continuous. (The displacement is not in general continuous on
  $\overline{\O(\Gamma)}^\infty$.)

In this section we describe a procedure for {\em calibrating} simplicial paths (see below
precise definitions).

\begin{defn}
	Let $\Sigma$ be a (simplicial) path in $\overline{\O(\Gamma)}^\infty$. We set $\lambda(\Sigma)$, the {\em displacement} of $\Sigma$, to be the supremum of displacements of points along $\Sigma$.
\end{defn}

\begin{defn}\label{defcal2}
  Let $L$ be a positive real number. A simplicial path $\Sigma=(X_i)_{i=0}^{k}$ in $\overline{\O(\Gamma)}^\infty$ is
  said to be $L$-{\em calibrated} if: 
  \begin{enumerate}[(i)]
  	\item $\lambda$ is continuous on $\Sigma$; 
        \item $\lambda(\Sigma)\leq L$;
  	\item no point $P$ of $\Sigma$ jumps (which, by Theorem~\ref{corlalx}, is equivalent to
          $\lambda(\phi)\leq\lambda(P)$);
        \item for any point $P$, in the interior of $\Sigma$ and that realises the maximum
          $\lambda(\Sigma)$, we have 
          $\lambda(P)=\lambda(\Delta_P)$ ({\em i.e.} $P$ is minimising in its simplex). Note
          that this implies that
          $\lambda(\Sigma)\in\operatorname{spec}(\phi)\cup\{\lambda(X_0),\lambda(X_k)\}$ (see
          Theorem~\ref{conj} for definition of $\operatorname{spec}(\phi)$). 
    \end{enumerate}
\end{defn}
\begin{rem}\label{remconv}
  If $A,B$ are two consecutive vertices of an  $L$-calibrated path then, by the continuity of $\lambda$, neither point can have jumped in the simplex they span. Hence by Lemma~\ref{lconv2} and property (ii) of Definition~\ref{defcal2},
  for any $P$ in the segment $\overline{AB}$ we 
  have, $$\lambda(\phi)\leq\lambda(P)\leq\max\{\lambda(A),\lambda(B)\}\leq L.$$ 
\end{rem}
\begin{thm}[Calibration]\label{thm3+10}
 Suppose Theorem~\ref{tconnected} is true in any rank less than $\rank(\Gamma)$.
	Let $\Sigma$ be a simplicial path in  $\overline{\O(\Gamma)}^\infty$ with finite
        displacement and 
        such that no point of $\Sigma$ jumps. Then in $\overline{\O(\Gamma)}^\infty$ there
        exists a $\lambda(\Sigma)$-calibrated simplicial path $\Sigma_o$ with
        the same endpoints as $\Sigma$.
\end{thm}
\proof

We outline the strategy of this proof to aid the reader.

\begin{itemize}
	\item First we regenerate $\Sigma$ to a path $\Sigma_1$ which lives inside
          $\O(\Gamma)$. This is basically an application Theorem~\ref{tregge} in its full
          generality, which, in particular, requires the inductive hypothesis on the rank. This
          is the only place of this section where such hypothesis is needed. We also note that
          in case $\phi$ is irreducible, then any path with finite displacement is in
          $\O(\Gamma)$ (see Remark~\ref{reminf}) so this step (and hence inductive hypothesis)
          is not  necessary in this case.
	\item Next, we define a simplicial path $\Sigma_2$ in $\overline{\O(\Gamma)}$, obtained from $\Sigma_1$ by, essentially, replacing each vertex with a point that minimizes the displacement in the corresponding simplex.
	\item Finally, we add extra points to $\Sigma_2$ in order to obtain a simplicial path $\Sigma_o$ to ensure that $\lambda$ is continuous along the path. 
	\item Along the way, we verify that we maintain control of the displacements of our paths, exploiting both quasi-convexity and the fact that $\operatorname{spec}(\phi)$ is well ordered (Theorem~\ref{conj}). 
\end{itemize}

Let $\Sigma=(X_i)_{i=1}^m$. Up to possibly adding extra vertices belonging to segments of
$\Sigma$, we may assume that it is alternating. (Note that this does not change the displacement of $\Sigma$).

Let $M=\min\{x\in\operatorname{spec}(\phi):\ x>\lambda(\Sigma)  \} $, which exists because
$\lambda(\Sigma)<+\infty$ and 
$\operatorname{spec}(\phi)$ is well ordered (Theorem~\ref{conj}). Let $\varepsilon>0$ so that $\lambda(\Sigma)+\varepsilon<M$.

We start by invoking Theorem~\ref{tregge}  to produce a
simplicial path $\Sigma_1=(W_j)_{j=1}^k$ in $\O(\Gamma)$, so
that $\lambda(\Sigma_1) \leq  \lambda(\Sigma) + \e   <M$ and so that $W_1$ and $W_k$ do not jump in in $\Delta_{X_1}$ and
$\Delta_{X_m}$ respectively. (Note that $\Delta_{X_1}$ is a face of $\Delta_{W_1}$, and $\Delta_{X_m}$
is of $\Delta_{W_k}$). 

We define a new simplicial path, $\Sigma_2$, as follows:
\begin{enumerate}
\item For any $j$, if $\Delta_{W_{j-1}}$ and $\Delta_{W_j}$ are both proper faces of some
  $\Delta_j$, then we add to  the path a new  point, $\widehat W_j\in\Delta_j\cap\overline{W_{j-1}W_j}$. We note that
  $\lambda(\widehat W_j) \leq \lambda(W_{j-1})$ and $\lambda(W_j)$, by quasi-convexity
  (Lemma~\ref{lconvexity}).

\item\label{stepp2}  We renumber the sequence of vertices, denoting them by $(W_j)_{j=1}^l$ (for some $l\geq
  k$). We now have a simplicial path which is alternating. 
\item\label{stepp3} For any $1\leq j\leq l$, we use Corollary~\ref{fatto2} and replace $W_j$ by a point $Y_j\in\overline{\Delta_{W_j}}^\infty$, chosen so that
  $\lambda(Y_j)=\lambda(\Delta_{W_j})=\lambda(\Delta_{Y_j})$, and requiring $Y_j$  to be
   maximal in the sense of Corollary~\ref{fatto2}.
\item\label{stepp4} We add endpoints $Y_0=X_1$ and $Y_{l+1}=X_m$.
\item\label{stepp5} If two consecutive points coincide, then we identify them and we
  renumber the sequence accordingly (and removing the corresponding segment). We call the
  resulting alternating simplicial path
  $\Sigma_2$. 
\end{enumerate}

\begin{lem}\label{lfine}
  For any vertex, $Y_j\in\Sigma_2$, we have that $\lambda(\Sigma)\geq\lambda(Y_j)\in\operatorname{spec}(\phi)\cup\{\lambda(X_1),\lambda(X_m)\}$.
\end{lem}
\proof The statement is obvious for endpoints. For other points, by construction, we have 
$M>  \lambda(\Sigma) + \e  \geq \lambda(W_j) \geq \lambda(Y_j)\in\operatorname{spec}(\phi)$, and our choice of $M$ implies $\lambda(Y_j)\leq\lambda(\Sigma)$. \qed

\begin{rem}\label{rembend}
 In Step~$(\ref{stepp4})$ we have $\lambda(Y_0)\geq \lambda(Y_1)$ and
 $\lambda(Y_{l+1})\geq\lambda(Y_l)$. This is because by definition $Y_1$ is the point in
 $\overline{\Delta_{X_1}}^\infty$ that realises
 $\lambda(\Delta_{X_1})=\inf_{T\in\Delta_{X_1}}\lambda(T)$, and $Y_0=X_1$. The same argument works
 for $Y_{l}$.
\end{rem}

\begin{lem}\label{lkk}
 Let $A,B$ be two consecutive vertices of $\Sigma_2$. Then, 
 \begin{enumerate}[(a)]
 \item For any point $P$ of $\overline{AB}$ we have
 $\lambda(P) \geq \lambda(\phi)$.
\item if  $\lambda(A)=\lambda(B)$, then
$\lambda$ is constant on the segment $\overline{AB}$;
\item if $\lambda(A)>\lambda(B)$, then
  there exists a simplex $\Delta\subset\O(\Gamma)$ and points $C,D$ so that:
  \begin{itemize}
  \item $A,B,C,D\in\overline{\Delta}^\infty$;
  \item $\lambda(A)<\lambda(C),\lambda(D)<\lambda(B)$;
  \item $\lambda$ is continuous on Euclidean segments $\overline{AC}$, $\overline{CD}$, and
    $\overline{DB}$. 
  \end{itemize}
\end{enumerate}
 \end{lem}
\proof Either $\lambda(A)\geq \lambda(B)$ or vice versa. Without loss of generality, up to
possibly switch the names of $A,B$, we may assume that $\lambda(A) \geq \lambda(B)$.

By how $\Sigma_2$ is defined, $A,B$ are introduced either in Step~$(\ref{stepp3})$ or
in Step~$(\ref{stepp4})$. Suppose first that both come from Step~$(\ref{stepp3})$. Since they
are consecutive  in $\Sigma_2$, by Step~$(\ref{stepp5})$ we may assume that the pair $\{A,B\}$
comes, in Step~$(\ref{stepp3})$, from a pair $\{W_j,W_{j+1}\}$ of two consecutive vertices of
the Step~$(\ref{stepp2})$-path. Since the path of Step~$(\ref{stepp2})$ is alternating, and
contained in $\O(\Gamma)$, either
$\Delta_{W_j}$ is a finitary face of $\Delta_{W_{j+1}}$ or vice versa: let $\Delta_0$ be the
one which is face of the other, and let $\Delta$ be the other. (We may have $\Delta=\Delta_0$.)
Moreover, from Step~$(\ref{stepp3})$, either $\lambda(A)=\lambda(\Delta_0)$ and
$\lambda(B)=\lambda(\Delta)$, or vice versa.

Since $\Delta_0$ is a face of $\Delta$ we have
$\lambda(\Delta_0)=\inf_{T\in\Delta_0}\lambda(T)\geq
\inf_{T\in\Delta}\lambda(T)=\lambda(\Delta)$. Since we are assuming $\lambda(A)\geq\lambda(B)$,
w.l.o.g. we may also assume $A\in\overline{\Delta_0}^\infty$ and
$B\in\overline{\Delta}^\infty$, thus $\lambda(A)=\lambda(\Delta_0), \lambda(B)=\lambda(\Delta)$.

Suppose now that $A$ is introduced in Step~$(\ref{stepp4})$. Then it is an endpoint, say
$A=Y_0=X_1$. In this case necessarily $B=Y_1$ is obtained in Step~$(\ref{stepp3})$ from
$X_1=A$. Whence $\lambda(A)>\lambda(\Delta_A)$ and $\lambda(B)=\lambda(\Delta_A)$. In this case
we set $\Delta_0=\Delta=\Delta_A$.

The same reasoning would work if $B$ is introduced in Step~$(\ref{stepp4})$, but then we would
get $\lambda(B)>\lambda(A)$ contradicting $\lambda(A)\geq\lambda(B)$. Therefore this latter
situation cannot happen. In particular $B$ is always introduced in Step~$(\ref{stepp3})$.

In any case, there there exists an open simplex $\Delta$ in $\O(\Gamma)$, with a (not
necessarily proper) finitary face, $\Delta_0$, such that $A \in \overline{\Delta_0}^\infty $, $B \in \overline{\Delta}^\infty$ and so that $\lambda(\Delta_0) \leq \lambda(A)$, $\lambda(\Delta) = \lambda(B)$. Thus both $A$ and $B$ belong to $\overline{\Delta}^\infty$.

Now let $\Delta_1$ be the simplicial face of $\Delta$ spanned by $A$ and $B$ (which may be
different from $\Delta$). Both $\lambda(A),\lambda(B)$ are finite. So, topologically, $A$ and $B$ are obtained from a graph, $X$, by collapsing invariant subgraphs $C_A$ and $C_B$, respectively. Therefore the points in $\Delta_1$ are obtained from $X$ by collapsing $C_A \cap C_B$, which is also invariant and hence all points in $\Delta_1$ have finite displacement. 

By the maximality of the dimension of $\Delta_B$    (Step $(\ref{stepp3})$),  and
Theorem~\ref{fatto2}, no point in $\Delta_1$ has jumped in $\Delta$. Hence, by
Theorem~\ref{newjump} and Lemma~\ref{lconv2}, for any point $P$, on the segment from $A$ to
$B$,
$$\lambda(\phi) \leq \lambda(\Delta) = \lambda(B) \leq \lambda(P) \leq \max \{ \lambda(A),
\lambda(B) \}.$$
This in particular proves $(a)$. Moreover, if $\lambda(A)=\lambda(B)$, we deduce that the
previous inequalities - except the first - are all equalities, thus proving $(b)$.

Finally, suppose that $\lambda(A) > \lambda(B)$. Since $\lambda$ is continuous in $\Delta_1$, and since $A$ has not jumped in $\Delta_1$ by Theorem~\ref{newjump}, we deduce - by Lemma~\ref{lconv2} -  that $\lambda$ is continuous along the segment from $A$ to $B$ except, possibly, at $B$.

If $\lambda$ is continuous in $\overline{AB}$, there is nothing to
prove. Otherwise, we use the fact that $B$ is defined in Step $(\ref{stepp3})$ by applying
Corollary~\ref{fatto2}. Our points $C,D$ correspond then to points $Z,W$ of
Corollary~\ref{fatto2}, which can be chosen with displacement arbitrarily close to
$\lambda(\Delta)=\lambda(B)$, in particular so that $\lambda(A)>\lambda(C),\lambda(D)$. The
fact that  $\lambda(C),\lambda(D)>\lambda(B)=\lambda(\Delta)$ follows from maximality condition
of $B$ (Step $(\ref{stepp3})$). Corollary~\ref{fatto2} also provides the continuity of $\lambda$ on the segments
$\overline{CD}$ and $\overline{DB}$. The continuity of $\overline{AC}$ follows from
Lemma~\ref{lconv2} because $A$ has higher displacement.

\qed

\medskip

We are now in position to finish the proof of Theorem~\ref{thm3+10}. Having $\Sigma_2$, we build
$\Sigma_o$ by using Lemma~\ref{lkk} to add points $C,D$ between consecutive vertices where
$\lambda$ is not continuous. 
In particular, $\lambda$ is continuous on $\Sigma_o$, and condition $(i)$ of
Definition~\ref{defcal2} is satisfied. Point $(a)$ of Lemma~\ref{lkk} gives condition $(iii)$. 

Note that added vertices are never point of maximum. Therefore Lemma~\ref{lfine}
provides condition $(iv)$. Finally Lemma~\ref{lfine} and Lemma~\ref{lconv2} imply that
$\lambda(\Sigma_o)\leq\lambda(\Sigma)$, so also condition $(ii)$ of Definition~\ref{defcal2} is
fulfilled with $L=\lambda(\Sigma)$. Thus $\Sigma_o$ is $\lambda(\Sigma)$-calibrated. 
\qed

\section{Preparation to peak reduction}
We keep Notation~\ref{not:gamma}. For the remaining of the section we fix $[\phi]\in\Out(\Gamma)$.
We recall that for simplices $\Delta\in\overline{\O(\Gamma)}^\infty$ we are using the notation $\lambda(\Delta)=\lambda_\phi(\Delta)=\inf_{X\in\Delta}\lambda_\phi(X)$.
In this section we prove some preliminary result needed to perform reduction of peaks.

We start by stating a (technical) fact that  can be informally
phrased as follows\footnote{We recall
  that by definition $\overline{\O(\Gamma)}^\infty=\overline{\O(\Gamma)}$ and that the symbol
  $\infty$ is just to put emphasis on the fact that we are considering the simplicial
  bordification of the outer space obtained by adding all simplices at infinity.}:

\centerline{\parbox{0.8\textwidth}{Given $X\in\overline{\O(\Gamma)}^\infty$ and $f:X\to X$
    an optimal map representing $[\phi]$, if $Y$ is sufficiently close to $X$ for the Euclidean
    metric, then any fold in $X$ directed by $f$ can be closely read in $Y$.}}

\begin{thm}\label{thmZnew}
Let $X,Y\in\overline{\Og(\Gamma)}$. Suppose that $\Delta_X$ is a simplicial face of
$\Delta_Y$. Thus as graphs, $X$ is obtained from $Y$ by collapsing a sub-graph $A$. Suppose that
$\core(A)$ is $\phi$-invariant. For $t\in[0,1]$ let $Y_t=(1-t)X +t Y$ be a
parametrization of the Euclidean segment from $X$ to $Y$. Let $\sigma_t:Y_t\to X$ be the map
obtained by collapsing  $A$ and by linearly rescaling the edges in $Y\setminus A$.

Let $f:X\to X$ be an optimal map representing $[\phi]$. Then for any $\e >0$ there is $t_\e>0$ such that $\forall 0\leq t<t_\e$ there is an optimal map
$g_t:Y_t\to Y_t$ representing $[\phi]$ such that $$d_\infty(\sigma_t\circ g_t,f\circ \sigma_t)<\e.$$
\end{thm}
\proof The proof of this theorem relies on accurate (but boring) estimates. For the happiness of
the reader we postpone the proof to the appendix.\qed

\begin{rem}
	Note that when $Y \in \overline{\O(\Gamma)}$, we may regard $\O(Y)$ as a subset of
        $\overline{\O(\Gamma)}$. Moreover, if $\lambda(Y) < \infty$, as is our usual
        assumption, then the same is true for all points in $\O(Y)$, since all points in this
        space share the same vertex groups which are necessarily invariant, by consequence of
        the fact that  $\lambda(Y) < \infty$. Note also that $\lambda$ is continuous on
        $\O(Y)$, because in general the displacement is continuous in the interior of {\em any}
        outer space. 
\end{rem}

\begin{rem}
	
	Consider the situation given by the hypotheses of \ref{thmZnew}. The $\phi$-invariance of $\core(A)$ allows us to build a straight map, $g:Y \to Y$, representing $[\phi]$ which leaves $\core(A)$ invariant. This map might not be optimal, but its Lipschitz constant provides an upper bound on the displacement of $Y$. 
	
	Now, along the path $Y_t$, we have the same topological trees (graphs of groups) except
        at the endpoint, $X$. We can thus re-scale edges but use the same topological straight
        map, $g$, to provide straight maps for all points $Y_t$ except for $X$. From the
        invariance of $\core(A)$, one easily sees that there is a constant, $C$, so that
        $\lambda(Y_t) < C$ for all points on the path. (We can include $X$ as well in this last statement).

	The hypotheses of Theorem~\ref{fatto1} therefore apply and we may deduce that $\lambda(X) \leq \liminf_{t \to 0} \lambda(Y_t) $.
\end{rem}

\begin{cor}\label{corpeak}
Let $X,Y\in\overline{\O(\Gamma)}$. We use the notation and hypotheses of
Theorem~\ref{thmZnew}. (In particular $\Delta_X$ is a simplicial face of $\Delta_Y$).
Let $f:X\to X$ be an optimal map representing $ [\phi]$. Suppose further that $\tau$ is an
$f$-illegal turn of $X$. Let $\Delta^\tau$ be the simplex obtained by folding $\tau$ and let
$X^\tau\in\Delta^\tau$ be the a point obtained from $X$ by folding $\tau$. 

\medskip

Given $\e > 0$, there exists $t_{\e}>0$, so that for all $t$ smaller than $t_{\e}$, there
exists an alternating simplicial path $\Sigma_t=(Z_i^t)_{i=0}^m$ in $\O(Y)$ from $Z_0^t=Y_{t}$ to a point $Z_m^t=Z_t$, so that 
\begin{itemize}
	\item $\Delta_{Z^t_i}$ has $\Delta_X$ as a simplicial face for any $i=0,\dots,m-1$,
        \item $\Delta_{Z_t}$ has $\Delta^\tau$ as a simplicial face,
	\item for any point $P$ of $\Sigma_t$ we have
          $\lambda(X)-\e<\lambda(P)\leq\lambda(Y_{t})$;
        \item for $s\in[0,t]$ the map $s\mapsto Z_s$ parametrizes the segment
          from $X^\tau$ to $Z_t$.
\end{itemize}

\end{cor}
\proof  For this proof we will work entirely with trees. So $Y$ will denote a $\Gamma$-forest,
$A$ an equivariant family of sub-trees --- that is to say, the full pre-image in $Y$ of an invariant
subgraph $\underline A\subseteq \underline Y$ --- and so on.

The map $\sigma_t$ is that introduced in the statement of Theorem~\ref{thmZnew}, and $g_t$ is the
map provided by Theorem~\ref{thmZnew}. Also, the $t_\epsilon$ is that provided by Theorem~\ref{thmZnew}.

We denote by $A_t$ the metric copy of $A$ in $Y_t$.
By hypothesis there are two different segments $\alpha_\tau,\beta_\tau$ incident at the same vertex $v$ in
$X$ such that $f$ overlaps $\alpha_\tau$ and $\beta_\tau$. If $v\notin\sigma_t(A_t)$ then, for
any small enough $\e$ and $t<t_\e$, also $g_t$ must overlap $\alpha=\sigma_t^{-1}(\alpha_\tau)$ and
$\beta=\sigma_t^{-1}(\beta_\tau)$, and the claim follows by (equivariantly) performing the
corresponding simple fold directed by $g_t$. Thus in this case the folding path consists of two
points: $Z_0^t=Y_t$ and $Z_1^t=Z_t$.
The inequality ``$\leq\lambda(Y_t)$'' follows
because the fold is directed by an optimal map, the inequality ``$>\lambda(X)-\e$'' follows by
lower semicontinuity of $\lambda$.

Otherwise, $\alpha$ and $\beta$ are segments incident to the
same component of $A_t$. If $\alpha$ and $\beta$ are incident to the same point, then we
proceed as above, so we can suppose that they are incident to different points of $A_t$.

For small enough $\e$ and $t<t_\e$ we
have that $g_t$ overlaps some open sub-segments of $\alpha$ and $\beta$. Let $a\in\alpha$
and $b\in\beta$ such that $g_t(a)=g_t(b)$ and such that $a$ is the closest possible to
$A_t$. Let $a'$ be the point where $\alpha$ is attached to $A_t$, and 
$b'$ the one where $\beta$ is attached to $A_t$.

Let $\gamma'$ be the segment from $a'$ to $b'$ in $A_t$, and let $\gamma$ be the segment between
$a$ and $b$ in $Y_t$. Clearly $\gamma=[a,a']\cup\gamma'\cup [b',b]$, where $[a,a']$ is the
sub-segment of $\alpha$ from $a$ to $a'$, and $[b',b]$ is the sub-segment of
$\beta$ from $b'$ to $b$. Note that $[a,a']\neq \alpha$ and $[b',b]\neq\beta$ because
$\alpha$ and $\beta$ are open and $a$ is the closest possible to $A_t$. 

On $\gamma$ we put an extra simplicial
structure given by the pull-back via $g_t$: we declare new vertices of $\gamma$ the points
whose $g_t$-image is a vertex of $Y_t$.
$g_t(\gamma)$ is a tree because $Y_t$ is. Moreover, since $g_t (a)=g_t(b)$,
the restriction of $g_t$ to $\gamma$ cannot be injective. In particular, if $x\in\gamma$ is a
point such that $d_{Y_t}(g_t(x),g_t(a))$ is maximal, then $x$ is a vertex of $\gamma$, and the two
sub-segments of $\gamma$ incident to $x$ are completely overlapped.

Let $Z_1^t$ be the tree obtained by equivariantly identify such segments. 
Note that $s\mapsto Z_1^s$ parametrizes the segment from $X$ to $Z_1^t$.
Clearly, $g_t$ induces
a map $g_t^1:Z_1^t\to Z_1^t$. Such map is continuous and not necessarily straight. However,
$$\Lip(g_t^1)\leq\Lip(g_t)$$ and $\PL(g_t^1)$ still represents $[\phi]$. Since
$\Lip(\PL(g_t^1))\leq\Lip(g_t^1)$ he have
$$\lambda(Z_1^t)\leq\lambda(Y_t).$$

Let $A'_t$ be the union of $A_t$ and the orbits of $[a,a']$ and $[b',b]$. Since
$[a,a']\neq\alpha$ and $[b',b]\neq \beta$, then the collapsing of $A'_t$ produces a point of $\Delta_X$.
As our identification occurred in $A_t'$, it follows 
that $\Delta_{Z_1^t}$ has $\Delta_X$ as a simplicial face.

Also, since $Y_t$ parametrizes the segment from $X$ to $Y$, as $t$ varies
$Z_1^t$ parametrizes the segment from $X$ to $Z_1^t$.

Note that a priori we may have $\Delta_{Z_1^t}=\Delta_Y$, but in any case
$\Delta_{Z_1^t}$ is either a (non necessarily proper)
simplicial face of $\Delta_Y$ or vice versa.

In $Z_1^t$ we have a simple path $\gamma_1$ resulting from $\gamma$ by the cancellation  of the
two identified segments at $x$. By construction $g_t^1$ is simplicial. If $g_t^1$ is not
injective on $\gamma_1$, we can iterate the above procedure and define points $Z_i^t$
with $$\lambda(Z_i^t)\leq \Lip(g_t)=\lambda(Y_t)$$
and such that $\Delta_{Z_i^t}$ has $\Delta_X$ as a simplicial face.
Moreover either
$\Delta_{Z_i^t}$ has $\Delta_{Z_{i-1}^t}$ as a simplicial face or vice versa, so the simplicial
path we are producing is alternating.
Since $\gamma$ has a finite number of vertices, we must stop, and we do when $\gamma_i$ is a
single point. At this stage, $\alpha$ and $\beta$ are incident to the same point and we are
reduced to the initial case. Note that any $Z_i^t\to X$ as $t\to 0$, thus so does any point
in segment from $Z_i^t$ to $Z_{i+1}^t$. Therefore by lower
semicontinuity of $\lambda$ for any $\e>0$, since we have
finitely many $Z_i^t$'s, for sufficiently small $t$ we have that for any $i$
$$\lambda(X)-\e<\lambda(Z_i^t)$$
and the same inequality holds for points in the segments from $Z_i^t$ to $Z_{i+1}^t$. Thus, up
to possibly replacing $t_\epsilon$ with a smaller positive number, we get that inequality of
third bullet in the statement, holds for any $t<t_\epsilon$.
\qed

\begin{rem}\label{newrem1}
  The length of the simplicial path produced by Corollary~\ref{corpeak} is bounded a priori by
  a constant depending only on $\rank(\Gamma)$. More precisely, consider the sequence of simplices
  $\Delta_{Z_i^t}$. It may happens that two consecutive $\Delta_{Z_i^t}$ and
  $\Delta_{Z_{i+1}^t}$ are equal, due to the fact that, in the proof of Corollary~\ref{corpeak},
  we subdivided  $\gamma$. Up to cancel such consecutive repetitions, the length of the
  sequence of $\Delta_{Z_i^t}$ is bounded by a constant depending on the complexity of $A_t$,
  hence on $\rank(\Gamma)$. 
\end{rem}

\begin{cor}\label{cordai}
	
Let $X,Y\in\overline{\O(\Gamma)}$ and suppose that $\Delta_X$ is a simplicial face of
$\Delta_Y$. Suppose that $\lambda(X) > \lambda(Y)$. 

Moreover, suppose that $X$ is an exit point for
$\Delta_X$\footnote{See Definition~\ref{exitp}}, and let $X_E$ be as
Definition~\ref{exitp}, chosen so that $\lambda(X_E) \geq \lambda(Y)$.  

Then there is a simplicial path $\Sigma=(W_i)$ in $\overline{\O(Y)}$, starting at $Y$ and
ending at $X_E$, with $W_i\in\O(Y)$ except possibly for the point $X_E$, such that 
for any point $P$ of $\Sigma$ we have 
$$\lambda(Y) \leq \lambda(P) \leq L< \lambda(X)$$
for some $L < \lambda(X)$.  
\end{cor}
\proof 
We inductively use Corollary~\ref{corpeak}: suppose that the exit point, $X_E$, is obtained by successive folds, $\tau_1, \ldots, \tau_m$.  (So that $\Delta_{X_E} = \Delta^{\tau_m}$.)

We parametrize the segment between $X$ and $Y$ by $Y_t=tY+(1-t)X$. 
 Lemma~\ref{lconvexity} and Lemma~\ref{lconv2} imply that on the Euclidean segment from $X$ to
 $Y$, the displacement is continuous, quasi-convex and strictly monotone near $X$. Hence, there
 exists a $t$ (which can be taken to be arbitrarily small), such that $Y_t$ satisfies
 $\lambda(X) - \e < \lambda(Y_t) < \lambda(X)$.  We then plug this in to
 Corollary~\ref{corpeak}, to find a point $Z_t$, whose displacement satisfies $\lambda(X) - \e <
 \lambda(Z_t) < \lambda(X)$, and a simplicial path, in $\O(Y)$, from $Y_t$ to $Z_t$, where all
 points met have the same displacement inequality,  where the path starts at $\Delta_Y$ and
 ends at $\Delta^{\tau_1}$. Since $s\mapsto Z_s$ parametrizes the segment from $X$ to $Z_t$, we
 are in position to apply Corollary~\ref{corpeak} again to the point $Z_t$, noting that
 $\Delta_X$ is a simplicial face of $\Delta_{Z_t}$ and that $\lambda(Z_t) < \lambda(X)$.  
 
 We continue inductively. 
 
  Concatenating our paths, and adding the points $Y$ and $X_E$, yields the result; the constant
  $L$ is simply the maximum displacement of points of our paths. By construction the
  displacement  is  a number strictly less than $\lambda(X)$ on vertices. Since $\Sigma\subset
  \O(Y)$ except possibly for its last point $X_E$, the displacement is continuous and
  quasi-convex (Lemma~\ref{lconvexity})  on $\Sigma$   except possibly at $X_E$ 
  where it may jump, but still lower-semicontinuity is preserved (Theorem~\ref{fatto1}). This
  implies that $L<\lambda(X)$. 

  \qed

  \begin{rem}\label{rempeakred}
As in Remark~\ref{newrem1}, up to repetitions, the simplicial length of the path $\Sigma$
provided by Corollary~\ref{cordai} is bounded a priori by a constant depending only on
$\rank(\Gamma)$. This is because of Remark~\ref{newrem1} and because the length of the path
from $X$ to $X_E$ is bounded by the dimension of $\O(\Gamma)$.    
  \end{rem}

\section{End of the proof of Theorem~\ref{tconnected}: peak reduction on simplicial paths}
We fix $\Gamma$ as in Notation~\ref{not:gamma} and $[\phi]\in\Out(\Gamma)$. Let
$\lambda=\lambda_\phi$.
We will prove:

\begin{lem}\label{lcsp}
  For any $L \geq \lambda(\phi) $, the level set
$$\{X\in\overline{\O(\Gamma)}^\infty:\lambda(\phi)\leq\lambda_\phi(X)\leq L \}$$
is connected by $L$-calibrated simplicial paths. 
\end{lem}
This in particular gives the second claim of Theorem~\ref{tconnected} (when $L=\lambda(\phi)$).
Moreover, if $\Sigma$  is any $L$-calibrated path  (hence in the above level set), then, by
possibly adding some extra vertices to $\Sigma$ we obtain a path in the same level set, and that  in
addition is alternating. So Theorem~\ref{tregge} applies and $\Sigma$ can be regenerated to
$\O(\Gamma)$, and this proves first claim of Theorem~\ref{tconnected}.

\medskip

We will proceed by induction and assume that Theorem~\ref{tconnected} is true in any rank less than $\rank(\Gamma)$.

\medskip

From now on we fix $A,B\in\overline{\O(\Gamma)}^\infty$ such that $\lambda(A),\lambda(B)\geq\lambda(\phi)$.
For any $L \geq \max\{\lambda(A),\lambda(B)\}$ we denote by $\Sigma_L(A,B)$ the set of $L$-calibrated simplicial paths from $A$ to $B$.

\begin{lem}
For some $L$, $\Sigma_L(A,B)\neq\emptyset$.
\end{lem}
\proof Since $\lambda(A),\lambda(B)\geq\lambda(\phi)$, they have not jumped.
Let $A'\in\Hor(A)$ and $B'\in\Hor(B)$, so that $A$ has not jumped in $\Delta_{A'}$ and $B$ has not jumped in $\Delta_{B'}$. Since $A',B'\in\O(\Gamma)$, which is connected,
there is a simplicial path in $\O(\Gamma)$ between $A',B'$. 
We can therefore use Theorem~\ref{thm3+10} to obtain
an element of $\Sigma_L$ (where the $L$ is the maximum displacement along such a path).\qed

\begin{defn}
  For any calibrated path $\Sigma=(X_i)$ we say that $X_i$ is a {\em peak} if
  $\lambda(X_i)=\lambda(\Sigma)$. A pair of two consecutive peaks $X_{i-1},X_{i}$ is called a {\em flat
    peak}. A peak is {\em strict} if it is not part of a flat peak.
\end{defn}

To any $\Sigma$ we can associate the triple $(\lambda(\Sigma),p,p_f)\in\operatorname{spec}(\phi)\times
\mathbb Z_{\geq 0}\times \mathbb Z_{\geq0}$ where $p$ is the number
of peaks, and $p_f$ that of flat peaks. We order $\operatorname{spec}(\phi)\times \mathbb Z_{\geq
  0}\times \mathbb Z_{\geq0}$ with lexicographic order, from left to right. That is, $(\lambda, p, p_f) > (\lambda', p', p_f')$ means:
\begin{itemize}
	\item $\lambda > \lambda'$, or
	\item $\lambda=\lambda'$ and $p > p'$, or 
	\item $\lambda=\lambda'$ and $p = p'$ and $p_f > p_f'$. 
\end{itemize}

\begin{lem}
There exists $\Sigma_0=(X_i)\in\Sigma_L(A,B)$, a calibrated path from $A$ to $B$, which
minimises $(\lambda,p,p_f)$. Namely, $\Sigma_0$ minimizes, in order:  
\begin{enumerate}
\item $\lambda(\Sigma)$;
\item the number peaks;
\item the number of flat peaks.
\end{enumerate}
\end{lem}
\proof By Theorem~\ref{conj} the set $\operatorname{spec}(\phi)$ is well-ordered, so
$\operatorname{spec}(\phi)\times \mathbb Z_{\geq 0}\times \mathbb Z_{\geq0}$ is
lexicographically well-ordered. Therefore every minimising sequence must eventually realise the
minimum.\qed

From now on we fix such a minimising $\Sigma_0$. 
\medskip

Note that if $X$ is a strict peak of a path $\Sigma$, then $\lambda$ is locally strictly
monotone near $X$, on both sides of $X$ in $\Sigma$. (By Lemma~\ref{lconvexity}.)

Once again, we need the inductive hypothesis.

\begin{lem}\label{newl2}
  Suppose that Theorem~\ref{tconnected} is true in any rank less than $\rank(\Gamma)$.
  Then $\Sigma_0$ has no strict peaks in its interior.
\end{lem}

\proof
Suppose that $\lambda(X_{i-1})<\lambda(X_i)>\lambda(X_{i+1})$. Set $X=X_i$, $Y=X_{i-1}, Z=X_{i+1}$, so that $\lambda(Y), \lambda(Z) < \lambda(X)$. 

By calibration, $X$ minimizes $\lambda$ in its simplex,
hence $\Delta_{X}$ is a proper face of both $\Delta_Y$
and $\Delta_Z$.

Since  $X$ is not a $\phi$-minimally displaced point, by Lemma~\ref{dom5}
$X\notin\TT(\phi)\subset\O(X)$. By
Lemma~\ref{LemmaX}, $X$ is an exit point. Let $X_E$ be as in Definition~\ref{exitp}. Since
$X_E$ can be chosen arbitrarily close to $X$, we  chose one so that $\lambda(X_E) \geq \max \{\lambda(Y), \lambda(Z)\}$.

Now we invoke Corollary~\ref{cordai} to get a simplicial path $\Sigma$ in $\overline{\O(Y)}$
from $Y$ to $X_E$, the displacement of whose points is between $\lambda(Y)$ and $L$, for some $L <
\lambda(X)$. In particular $\lambda(\Sigma)<\lambda(X)$. 

We now interpret this as a simplicial path in $\overline{\O(\Gamma)}$.
Since $\lambda(Y) \geq \lambda(\phi)$ no point of such path jumps. We apply
Theorem~\ref{thm3+10} to obtain a calibrated path $\Sigma_Y$ from $Y$ to $X_E$, whose
displacement is less than 
 $\lambda(X)$. By symmetry, we get a calibrated path $\Sigma_Z$ from $X_E$ to $Z$ whose displacement is
 less than $\lambda(X)$. Let $\Sigma_1$ be the simplicial path obtained by following $\Sigma_0$
 till $Y$, then $\Sigma_Y$, then $\Sigma_Z$ and then again $\Sigma_0$ till its end. Since
 $\lambda(\Sigma_Y),\lambda(\Sigma_Z)<\lambda(X)=\lambda(\Sigma_0)$, we have 
  $\lambda(\Sigma_1)\leq\lambda(\Sigma_0)$.

If $\lambda(\Sigma_1)<\lambda(\Sigma_0)$, we apply Theorem~\ref{thm3+10} and contradict the
minimality of $\Sigma_0$. Otherwise, paths $\Sigma_Y$ and $\Sigma_Z$ do not contain peaks
of $\Sigma_1$. Therefore $\Sigma_1$ is a $\lambda(\Sigma_0)$-calibrated which has fewer strict peaks
than $\Sigma_0$, contradicting minimality.
\qed

\begin{lem}\label{newl1}
  $\Sigma_0$ has no flat peaks unless $\lambda$ is constant on $\Sigma_0$ and $\lambda(\Sigma_0)=\lambda(\phi)$.
\end{lem}
\proof If the function $\lambda$ is not constantly equal to $\lambda(\phi)$ on $\Sigma_0$, then in particular $\lambda$ is
strictly bigger than $\lambda(\phi)$ on peaks.
Suppose that there is $Y,X$ two consecutive vertices of $\Sigma_0$ with
$$\lambda(X)=\lambda(Y)=\lambda(\Sigma_0)>\lambda(\phi).$$ The idea is to find a third point $Z$ to
add between $Y$ and $X$ in order to destroy the flat peak.
If there is a point $Z$ in the interior of the segment $YX$, with $\lambda(\phi)\leq\lambda(Z)<\lambda(X)=\lambda(Y)$, then we
add it.

Otherwise, $\lambda$ is constant on $\overline{XY}$. Let $W$ be a point in the interior of the
segment $\overline{XY}$. If $W$ is not a local minimum for $\lambda$ in $\Delta_W$, then near
$W$ we find $Z$ with the above properties. We add it.

If $W$ is a local minimum for $\lambda$ in $\Delta_W$ then, by Lemma~\ref{dom5} and
Lemma~\ref{LemmaX}, near $W$ in $\O(W)$  there is a point $Z$ with the above properties and such that $\Delta_W$
is a finitary face of $\Delta_Z$ in $\O(W)$. We add $Z$.

In each case, we have added a point, $Z$, such that $\Delta_X$ and $\Delta_Y$ are faces of $\Delta_Z$, and since the original path was calibrated, we can verify - using Theorem~\ref{newjump} - in each case that $X,Y$ did not jump in $\Delta_Z$. Hence we can add $Z$ to the path. By Lemma~\ref{lconv2}, the new path is still a calibrated path (continuity at $Z$ is automatic, since $\lambda$ is continuous in $\O(W)$), with the same displacement as $\Sigma_0$, and the same number of peaks, but with one less flat peak, contradicting the minimality of $\Sigma_0$.
\qed

\medskip
It follows that the maximum displacement of points of $\Sigma_0$ is  reached at endpoints. Thus
$\Sigma_0$ is a calibrated simplicial path in the requested level set, proving Lemma~\ref{lcsp}.
To finish the proof of Theorem~\ref{tconnected}, simply observe that we have shown that we can
connect any two points in $\{X\in\overline{\O(\Gamma)}^\infty:\lambda_\phi(X)=\lambda(\phi)\}$
by a calibrated simplicial path with no peaks, either strict or  flat, unless the displacement
is constant. This immediately implies that the displacement is constant along the path. \qed

\medskip

The following is an observation that may be helpful for algorithmic purposes.

\begin{rem} If $\phi$ is irreducible,
  there exists a constant $K$, depending only on $\rank(\Gamma)$, such that,
  given a $L$-calibrated alternating simplicial path $\Sigma$ having some peak in its interior,
  and such 
  that the displacement is not constant along $\Sigma$, there exists a $L$-calibrated alternating
  simplicial path $\Sigma'$ with either less displacement or one peak less, and whose 
  simplicial length is increased at  
  most by $K$.

  This is because  we can remove a strict peak from $\Sigma$ as in Lemma~\ref{newl2}
  --- if $\Sigma$   contains no strict peak, we create one as in
  Lemma~\ref{newl1}, without changing the global number of peaks nor $\lambda(\Sigma)$, and
  increasing the length of 
  $\Sigma$ by $1$ ---. The
  control on simplicial length comes from the use of Corollary~\ref{cordai} and Theorem~\ref{thm3+10} in
  the proof of Lemma~\ref{newl2}:

  By
  Remark~\ref{rempeakred} any use of Corollary~\ref{cordai} increase the simplicial length
  by a fixed amount, and  since $\phi$ is irreducible, every calibrated path is in
  $O(\Gamma)$; therefore the calibration process Theorem~\ref{thm3+10} does not involve
  regeneration of paths, nor continuity issues, (so the alternating $\Sigma_2$ is already
  calibrated in the proof of Theorem~\ref{thm3+10}), and it is readily checked that in this
  case  calibration increases the  length by a fixed amount.   
\end{rem}

\section{Applications}\label{appl}

In this section we show how the connectedness of the level sets gives a solution to some
decision problems. Namely we will prove Theorems~\ref{detectirred} and~\ref{conjirred} and some
generalisations.
We will work with {\bf graphs in the volume-one slice of $CV_n$.}

Recall that a point, $X$, of $CV_n$ is called $\e$-thin
if there is a homotopically non-trivial loop in $X$ of length at most $\e$. Conversely, $X$ is
called $\e$-thick if it is not $\e$-thin.

\begin{prop}[{\cite[Proposition 10]{BestvinaBers}. See also~\cite[Proposition 5.5]{FMpartI},
    and~\cite[Section 8]{FM13}}] \label{thick}
Let $X \in CV_n$ (that is, $X$ is a volume-one marked metric graph) and $f:X \to X$ a straight
map representing
some automorphism of $F_n$. Let $\lambda  = Lip(f)$, let $N$ equal the maximal length of chains
of topological subgraphs of any graph in $CV_n$ (this is clearly a finite number) and let $\mu$
be any real number greater than $\lambda$. Then if $X$ is $1/((3n-3)\mu^{(N+1)})$-thin, then it
has a nontrivial core sub-graph which is $f$-invariant up to homotopy, in particular the automorphism represented by $f$ is reducible. For instance, one can take $N=3n-3$.
\end{prop}

\begin{defn}
A uniform rose in $CV_n$ is a rose-graph ({\em i.e.} a bouquet of circles) whose edges all have the same length. Let $X \in CV_n$. Then
we call $R$ an adjacent uniform rose if it obtained by collapsing a maximal tree in $X$ and
then rescaling so that all edges in $R$ have the same length.
\end{defn}

\begin{prop}
Let $X \in CV_n$ be a point which is $\e$-thick and let $R$ be any adjacent uniform rose (both of volume 1). 
Then, $\Lambda(X, R) \leq 1 / \e$ and $\Lambda(R, X) \leq n$. 
\end{prop}
\proof By Theorem~\ref{sausagelemma}, we can look at candidates that realise the stretching
factor. Since, topologically, one passes from $X$ to $R$ by collapsing a maximal tree, we get
that a candidate in $X$, when mapped to $R$, crosses every edge at most twice. In fact the
candidate crosses every edge of $R$ at most once in the case of an embedded simple loop or an
infinity-symbol loop. This gives the first inequality, on taking into account that $X$ is $\e$-thick and that barbells have length at least $2 \e$. 

For the second inequality note that an embedded loop in $R$ is an edge and has length $1/n$ and
lifts to an embedded loop in $X$, of length at most $1$. An infinity-symbol loop in $R$ consists of two distinct edges, has length $2/n$ and lifts to a loop in $X$ which goes through every edge at most twice. (Barbells are not present in $R$). 
\qed

\begin{cor} \label{roseapprox}
Let $X \in CV_n$ be $\e$-thick and let $R$ be an adjacent uniform rose. Consider $[\phi] \in \Out(F_n)$. Then $\Lambda(R, \phi R) \leq \frac{ n }{\e}\Lambda(X, \phi X)$.  
\end{cor}

Now, we use connectedness of level sets (Theorem~\ref{tconnected}) for deducing the following result.
\begin{prop}\label{roseconnect}
Let $R, R_{\infty}$ be two points in $CV_n$ which are both uniform roses. Let $[\phi] \in \Out(F_n)$ be irreducible and suppose that $\mu$ is any real number greater than: 

$$\max\{ \Lambda(R, \phi R), \Lambda(R_{\infty}, \phi R_{\infty})\}.$$

Then there exist $R_0=R, R_1, R_2, \ldots, R_k=R_{\infty}$, which are all uniform roses in $CV_n$ such that: 
\begin{itemize}
\item For each $i$, there exists a simplex $\Delta_i$ such that $\Delta_{R_i}$ is a rose-face of both  $\Delta_i$ and $\Delta_{i+1}$. 
\item $\Lambda(R_i, \phi R_i) \leq \frac{ n }{\e} \mu$, where $\e = 1/((3n-3)\mu^{(N+1)})$. 
\end{itemize}
\end{prop}
\proof
This follows from Theorem~\ref{tconnected}, using Definition~\ref{def:sp}, since each pair $\Delta_i$ and $\Delta_{i+1}$ have a (at least one) common rose face; just take any uniform adjacent rose in any  common rose face. The remaining point follows from Corollary~\ref{roseapprox} and Proposition~\ref{thick}. 
\qed

\ \medskip

\noindent
{\em Proof of Theorem~\ref{conjirred}}: We clearly have an algorithm which terminates (Remark~\ref{remS}), and it is apparent that if $\psi \in S_{\phi}$ then these automorphisms are conjugate. It remains to show the converse; that if they are conjugate, then $\psi \in S_{\phi}$. 

Let $R$ be the uniform rose corresponding to the basis $B$. If $\psi$ were conjugate to $\phi$,
then there would be a conjugator, some $[\tau] \in \Out(F_n)$ such that $\psi = \tau \phi
\tau^{-1}$. Let $R_{\infty} = \tau R$. Remind that the $\Out(F_n)$-action on $CV_n$ is a
right-action, namely $\phi(\psi(X))=(\psi\phi)X$. In
particular,
\begin{equation}\label{eqb}
||\psi||_B=\Lambda(R,\psi
R)=\Lambda(R,(\tau\phi\tau^{-1})R)=\Lambda(\tau^{-1}(\tau R),\tau^{-1}(\phi(\tau
R)))=\Lambda(R_\infty,\phi R_\infty). 
\end{equation}

Now we use Proposition~\ref{roseconnect} to find a sequence $R=R_0, R_1, \ldots, R_k=R_{\infty}$, such that each consecutive pair are incident to a common simplex and $\Lambda(R_i, \phi R_i) \leq n(3n-3) \mu^{3n-1} = K$. 

Let $\tau_i$ so that $R_i=\tau_iR$.  Since $R_i$ and $R_{i+1}$ are both
incident to a common simplex, there exists a CMT automorphism $\zeta_i$ such that
$\tau_i(\zeta_i(\tau_i^{-1}( R_i)))=R_{i+1}$.
Thus $$\zeta_i\tau_iR=\tau_i(\zeta_i(R))=\tau_i(\zeta_i(\tau_i^{-1}(R_i)))=R_{i+1}=\tau_{i+1}R,$$
and up possibly
compose $\zeta_i$ with a graph-automorphism of $R$, we may assume
$\tau_{i+1}=\zeta_i\tau_i$. Therefore $\tau_{i+1}=\zeta_i\dots\zeta_0$ (and we set $\tau_0=id$).

Now let $\phi_i =  \tau_i \phi \tau_i^{-1}$. Clearly $\phi_0=\phi$ and $\phi_k=\psi$.

Since $\phi_{i+1} = \zeta_i \phi_i \zeta_i^{-1}$, to finish the proof we just need that
$||\phi_i ||_B \leq K$. This follows since, as in \eqref{eqb} 
$$
||\phi_i||_B= \Lambda(R, \phi_i R) = \Lambda( R_i, \phi R_i ) \leq K.
$$

\qed

We prove now Theorem~\ref{detectirred}. First a lemma, 

\begin{lem}\label{reducerose}
Let $X$ be a core graph and $f$ a homotopy equivalence on $X$, having a proper subgraph $X_0$,
with nontrivial fundamental group, such that $f(X_0) = X_0$.  Then there is a maximal tree, $T$, such that the automorphism induced by $f$ on the rose $X/T$ is visibly reducible.  
\end{lem}
\proof
Choose $X_0$ to be minimal. Therefore it will have components, $X_1, \ldots, X_k$ such that $f(X_i)=X_{i+1}$ with subscripts taken modulo $k$. Take a maximal tree for each $X_i$ and extend this to a maximal tree, $T$, for $X$. It is then clear that if we take $B_i$ to be the set of edges in $X/T$ coming from $X_i$, that $f_*$ will be visibly reducible as witnessed by $B_1, \ldots, B_k$. (Note each subgroups generated by each $B_i$ are only permuted/preserved up to conjugacy, since the $X_i$ are disjoint and so one cannot choose a common basepoint).  
\qed

\medskip

\noindent
{\em Proof of Theorem~\ref{detectirred}}: 
The algorithm clearly terminates (Remark~\ref{remS}), and if there is a $\psi$ in $S^+$ which is visibly reducible, then $\phi$ is reducible. It remains, therefore, to show that if $\phi$ is reducible, then there is some $\psi \in S^+$ which is visibly reducible.

We proceed much as in the proof of Theorem~\ref{conjirred}, but here we do not know that the
points in $CV_n$ we encounter will remain uniformly thick.

Let $R$ be the uniform rose corresponding to the basis $B$. By Theorem~\ref{strongcorred},
there exists an $X \in CV_n$ with a core invariant subgraph and such that $\Lambda(X, \phi(X))
< \mu$. 
By  Theorem~\ref{tconnected}, there exists a simplicial path from $R$ to $X$, whose vertices are   
points, $X_0=R, X_1, \ldots, X_k=X$, such that $\Lambda(X, \phi(X_i)) < \mu$. Choose the
maximal index, $M$, such that $X_0, X_1, \ldots, X_M$ are all $\e$-thick, where $\e =
1/((3n-3)\mu^{(N+1)})$ as in Proposition~\ref{thick}.

If $M<k$, then $X_{M+1}$ is $\varepsilon$-thin, and by Proposition~\ref{thick}, we have that
$X_{M+1}$ has an optimal representative for $[\phi]$ which admits an invariant
subgraph. Therefore, up to replacing $X$ with $X_{M+1}$, we may assume that $X_i$ is
$\varepsilon$-thick for $i=0,\dots,k-1$.

Since $X_k$ has an invariant subgraph,  by
Lemma~\ref{reducerose}, we may find an adjacent uniform rose face, $R_k$, so that the
representative of $[\phi]$ at $R_k$ is visibly reducible.

Now, for each $i \leq k-1$, we find a uniform rose $R_i$ which
is adjacent to both $X_i$ and $X_{i+1}$, which exist by definition of simplicial
path (Definition~\ref{def:sp}). Note that since $X_0=R$ is a rose, then $R_0=R$. Moreover, by
Corollary~\ref{roseapprox} we have $\Lambda(R_i,\phi R_i)<K$ for any $i=0,\dots, k-1$.

We now conclude exactly as in the proof of Theorem~\ref{conjirred}: Let $[\tau] \in \Out(F_n)$ be such that $R_k=\tau R$, and let
$\psi = \tau \phi \tau^{-1}$. Find CMT automorphisms $\zeta_i$ such that
$\tau_i=\zeta_{i-1}\dots\zeta_0$ satisfies $R_i=\tau_iR$ and $\tau_{k}=\tau$. Define $\tau_0=id$ and  $\phi_i =  \tau_i \phi \tau_i^{-1}$, so that $\phi_0=\phi$,
$\phi_{k} = \psi$, and $\phi_{i+1}=\zeta_i\phi_i\zeta_i^{-1}$.

Since each $\Lambda(R_i,\phi R_i)<K$, as in~\eqref{eqb}, we get  that each $\phi_i \in S_i$ for $i
\leq k-1$. Hence $\psi \in S^+$ and is visibly reducible, as desired.
\qed

\subsection{Generalisations}\label{s2.1}
Our algorithms work in some more general setting that just free groups. For instance, consider
the case of a group $G$ 
equipped with a splitting $\G=(\{G_i\},n)$ where the factor groups $G_i$ are finite groups.
In this case $\Og(\G)$ is a deformation space of finite graphs of groups with trivial
edge-groups and finite vertex groups.

This leads to Theorem~\ref{finite}, which we now explain how to prove.

Theorems~\ref{conjirred} and~\ref{detectirred} generalise as follows. As above, we work in the
volume-one slice of $\Og(\G)$.  Instead of uniform roses
one can use uniform ``hairy roses'', that is to say, graph $X\in\Og(\G)$ obtained from a rose by
attaching, to the unique vertex, edges each ending with a non-free vertex. Uniform here means that all edges have the same
length. 

Any $X\in\Og(\G)$ is a face of a simplex
containing a hairy rose simplex: to see this, first, for any non-free vertex $v$ which is not a
leaf, fold a little all edges at $v$; then, once all non-free vertices are leaves, collapse a
maximal tree in the sub-graph consisting of edges incident only at free vertices.
We say that a uniform hairy rose is adjacent to $X$ if obtained in this way, plus a rescaling
of edges. 

Now define a CMT automorphisms of a hairy rose as a change of marking between two hairy roses
`adjacent' to a common point.
More precisely, let $\Delta_1,\Delta_2$ be simplices in $\Og(\G)$, having a common face; let
$X\in\overline{\Delta_1}\cap\overline{\Delta_2}$,  and let $R_1,R_2$ be uniform hairy roses in
$\overline{\Delta_1},\overline{\Delta_2}$ respectively. Then we call $R_1$ and $R_2$ adjacent. 

Then, letting $R$ be a fixed marked hairy rose, we define,	
$$
CMT_R(G) = \{ [\phi] \in \Out(G) \ :  \ \phi(R) \text{ is adjacent to } R \}.
$$ 

\begin{rem}
	We note that this slightly different to the notion of adjacency in $CV_n$, but the idea is very similar. We start with an alternating simplicial path and want to replace each vertex along that path with a hairy rose. In $CV_n$, one can do this by replacing each point with a rose in such a way that consecutive roses are in faces of a common simplex. In this situation, moving to a hairy rose involves inserting `stems' and then collapsing a maximal tree (ignoring the stems). However, there are several (although finitely many) ways of introducing these stems since the vertex groups are non-trivial. This means each vertex in the original simplicial path gives rise to two hairy roses - one can insert stems consistently between consecutive points, but not necessarily for three consecutive points - and in the resulting sequence of hairy roses, consecutive hairy roses are adjacent in the sense described above. 
\end{rem}


There are finitely many CMT automorphisms since the finiteness of the vertex groups implies that the stabiliser of any point is finite, and also that the deformation space  $\Og(\G)$ is locally finite (and so there are only finitely many hairy roses adjacent to a given one). Moreover, since $\Og(\G)$ is
connected, the CMT automorphisms generate $\Out(\G)$.

Now we can build algorithms exactly as in Theorems~\ref{conjirred} and~\ref{detectirred}. The
fact that vertex groups are finite implies that Remark~\ref{remfin} holds true. So the set $S$
in the statements is finite, and algorithms stop in finite time. The fact that there are finitely many 
CMT automorphisms implies that the set $S^+$ in Theorem~\ref{detectirred} is finite.

The proof that these algorithms work now goes {\em mutatis mutandis} as in the case of $CV_n$.
In particular, the conjugacy problem for irreducible automorphisms and the detection of reducibility are solvable in $\Out(\G)$. 
\newpage
\section{Appendix: proof of Theorem~\ref{thmZnew}}\label{appendix}

In this section we give the proof of Theorem~\ref{thmZnew}, which we restate for convenience (recall we are
using Notation~\ref{not:gamma} and $[\phi]\in\Out(\Gamma)$).
\begin{thm*}[Theorem~\ref{thmZnew}]
Let $X,Y\in\overline{\Og(\Gamma)}$. Suppose that $\Delta_X$ is a simplicial face of
$\Delta_Y$. Thus as graphs, $Y$ is obtained by collapsing a sub-graph $A$. Suppose that
$\core(A)$ is $\phi$-invariant. For $t\in[0,1]$ let $Y_t=(1-t)X +t Y$ be a
parametrization of the Euclidean segment from $X$ to $Y$. Let $\sigma_t:Y_t\to X$ be the map
obtained by collapsing  $A$ and by linearly rescaling the edges in $Y\setminus A$.

Let $f:X\to X$ be an optimal map representing $[\phi]$. Then for any $\e >0$ there is $t_\e>0$ such that $\forall 0\leq t<t_\e$ there is an optimal map
$g_t:Y_t\to Y_t$ representing $[\phi]$ such that $$d_\infty(\sigma_t\circ g_t,f\circ \sigma_t)<\e.$$
\end{thm*}
\proof We split the proof in two sub-cases. First when $A$ is itself a core graph, and then the
case when $\core(A)$ is trivial. Clearly the disjoint union of the two cases implies the mixed case.

We will work at once with graphs and trees, by using Notation~\ref{tildeunderbar}.

\begin{lem}[When $A$ is a core graph]\label{lemmaZnew}
Let $X,Y\in\overline{\Og(\Gamma)}$. Suppose that as graphs of groups, $X$ is obtained from $Y$
by collapsing a $\phi$-invariant core sub-graph $A=\sqcup A_i$.
Then the conclusion of Theorem~\ref{thmZnew} holds.
\end{lem}
\proof  We begin by fixing  some notation. First of all, we will use the symbol $\lambda$ to
denote any of the displacement functions of $\phi$ (i.e. $\lambda_\phi,\lambda_{\phi|_A},\dots$).
If $x$ is a point in a metric space, we denote by $B_r(x)$ the open metric ball centered at $x$
and radius $r$.
For any $i$, we denote by $v_i$ the non-free vertex of $X$
obtained by collapsing $A_i$. For any $t$ we denote by $A^t$ the metric copy of $A$ in
$Y_t$. Note that $A$ is uniformly collapsed in $Y_t$, that is to say, $[A^t]\in\mathbb P\O(A)$
is the same element for any $0<t\leq 1$, and we have $\vol(A^t)=t\vol(A^1)$.

By lower
semicontinuity of $\lambda$ (Theorem~\ref{fatto1}) we have that
\begin{equation}
  \label{eq:0}
  \forall \e_0>0\exists t_{\e_0}>0 \text{ such that } \forall t<t_{\e_0} \text{ we have }
\lambda (Y_t)>\frac{\lambda (X)}{1+\e_0}.
\end{equation}

A priori $f$ may collapse some edge, in any case   $\forall \e_1>0\exists f_1:X\to X$
a straight map representing $[\phi]$ such that $f_1$ does not collapse any edge, and
\begin{equation}
  \label{eq:1}
  d_\infty(f,f_1)<\e_1\quad \text{and}\quad \Lip(f_1)<\Lip(f)(1+\e_1)=\lambda(X)(1+\e_1).
\end{equation}

Moreover $\exists 0<\rho_0=\rho_0(X,f_1)$ such that $\forall \rho<\rho_0$
\begin{itemize}
\item $B_\rho(x)$ is star-shaped for any $x\in X$ (i.e. it contains at most one vertex);
\item for any $i$, each connected component of $f_1^{-1}(B_\rho(v_i))$ is star-shaped
  and contains exactly one pre-image of $v_i$;
\item for any $i,j$ the connected components of
  $f_1^{-1}(B_\rho(v_i))$ and those of $f_1^{-1}(B_\rho(v_j))$ are pairwise disjoint.
\end{itemize}

We fix an optimal map $\varphi:A^1\to A^1$ representing $[\phi|_A]$. Since $[A^t]\in\mathbb
P\O(A)$ does not depend on $t$, $\varphi:A^t\to A^t$ is an optimal map for any $t\in(0,1]$ and the
Lipschitz constant does not change. Clearly (by Sausage
Lemma~\ref{sausagelemma})
\begin{equation}
  \label{eq:fi}
   \Lip(\varphi)\leq \lambda(Y_t)\quad \text{ for any } t.
\end{equation}
The natural option is to define $g_t$ by using $\sigma_t^{-1}\circ
f_1\circ\sigma_t$. Hence, we need to deal with places where $\sigma_t^{-1}$ is not
defined. (We have to understand how to deal with arcs in $X$ whose $f_1$-image crosses
some $v_i$.)

We fix lifts $\wt\varphi$ of $\varphi$ and $\wt f_1$ of $f_1$.
For any $v_i$, and any $x\in f_1^{-1}(v_i)$, to any germ of edge $\alpha$ at $x$  we associate
a path $\gamma_\alpha\in Y$ as follows. We do two different constructions: one in case $x$
is one of the $v_j$'s, and another for the case when $x$ is different from others $v_j$'s.

{\bf Case 1}. Suppose $x=v_k$ and $f_1(x)=v_i$ for some $k,i$ (not necessarily different).
Let $\alpha$ be a germ of edge at $x$. First of all we {\bf choose} a lift
$\wt\alpha$ of $\alpha$. All subsequent choices of lifts of objects, made during the definition of
$\gamma_\alpha$, will depend on, and will be uniquely determined by, the choice of
$\wt\alpha$. After having defined $\gamma_\alpha$, we forget about such choices of lifts.

The  germ $\alpha$ corresponds to a germ $\alpha_Y(=\sigma^{-1}_t(\alpha))$ in $Y$ incident to
$A_k$ at a point that we denote by $p_\alpha$.  The lift $\wt\alpha$ corresponds to a germ
$\wt\alpha_Y$ incident to $\wt p_\alpha\in\wt A_k$, where $\wt p_\alpha$ is a preimage of $p_\alpha$ and $\wt A_k$ is the component of the
preimage of $A_k$ containing $\wt p_\alpha$. (See Figure~\ref{pippo}.)

\setlength{\unitlength}{1ex}
\begin{figure}[htbp]
  \centering
  \begin{picture}(70,38)
\put(8,0){
    \put(5,5){\circle{10}}
    \put(5,1.1){\line(0,1){7.8}}
    \put(8.9,5){\line(1,0){10}}
    \put(8.9,5){\makebox(0,0){$\bullet$}}
    \put(2,0){\makebox(0,0){$A_k$}}
    \put(10,3){\makebox(0,0){$p_\alpha$}}
    \put(15,6){\makebox(0,0){$\alpha_Y$}}

    \put(40,0){
      \put(-0.25,5){\circle{5.85}}
      \put(5.9,5){\circle{5.85}}
      \put(8.9,5){\line(1,0){10}}
      \put(8.9,5){\makebox(0,0){$\bullet$}}
      \put(2,0){\makebox(0,0){$A_i$}}
      \put(11,3){\makebox(0,0){$p_\beta$}}
      \put(15,6.5){\makebox(0,0){$\beta_Y$}}
      \put(-2,0){\multiput(0,0)(6,0){3}{
                  \put(-1,14){\line(1,1){5}}
                  \put(1,20){\line(1,-1){6}}}
                \put(16.5,16.5){\makebox(0,0){$\bullet$}}
                \put(16.5,16.5){\line(1,0){5}}
                \put(15.5,14.5){\makebox(0,0){$\wt p_\beta$}}
                \put(19,18){\makebox(0,0){$\wt \beta_Y$}}
                \put(1,16){\makebox(0,0){$\bullet$}}
                \put(-2,17){\makebox(0,0){$\wt\varphi(\wt p_\alpha)$}}
      \put(2,12){\makebox(0,0){$\wt A_i$}}
                }
    \put(5,25){\line(5,1){14}}
    \put(5,25){\makebox(0,0){$\bullet$}}
    \put(11,29){\makebox(0,0){$\wt \beta=\wt f_1(\wt \alpha)$}}
    \put(2,25){\makebox(0,0){$\wt v_i$}}

    \put(0,8){
      \put(5,25){\line(5,1){14}}
      \put(5,25){\makebox(0,0){$\bullet$}}
      \put(11,29){\makebox(0,0){$\beta=f_1(\alpha)$}}
      \put(2,25){\makebox(0,0){$v_i$}}
               }

      }

    \multiput(0,0)(6,0){3}{
       \put(-1,14){\line(1,1){4}}
       \put(1,20){\line(1,-1){5}}}
      \put(16.5,16.5){\makebox(0,0){$\bullet$}}
      \put(16.5,16.5){\line(1,0){5}}
      \put(15.5,14.5){\makebox(0,0){$\wt p_\alpha$}}
      \put(19,18){\makebox(0,0){$\wt \alpha_Y$}}
      \put(2,12){\makebox(0,0){$\wt A_k$}}

    \put(5,25){\line(5,1){10}}
    \put(5,25){\makebox(0,0){$\bullet$}}
    \put(10,28){\makebox(0,0){$\wt\alpha$}}
    \put(5,28){\makebox(0,0){$\wt x$}}

    \put(0,8){
      \put(5,25){\line(5,1){10}}
      \put(5,25){\makebox(0,0){$\bullet$}}
      \put(10,28){\makebox(0,0){$\alpha$}}
      \put(5,28){\makebox(0,0){$x$}}
               }
\put(30,25){\makebox(0,0){$\wt\gamma_\alpha$}}
\put(31.5,24){\vector(2,-1){14}}
}
\put(0,5){in $Y$:}
\put(0,15){in $\wt Y$:}
\put(0,25){in $\wt X$:}
\put(0,33){in $X$:}
\put(33,30){\vector(1,0){6}}
\put(36,32){\makebox(0,0){$f_1$}}
\put(33,10){\vector(1,0){6}}
\put(36,12){\makebox(0,0){$\varphi$}}

  \end{picture}
  \caption{How to  choose the paths $\wt\gamma_\alpha$}\label{pippo}
\end{figure}
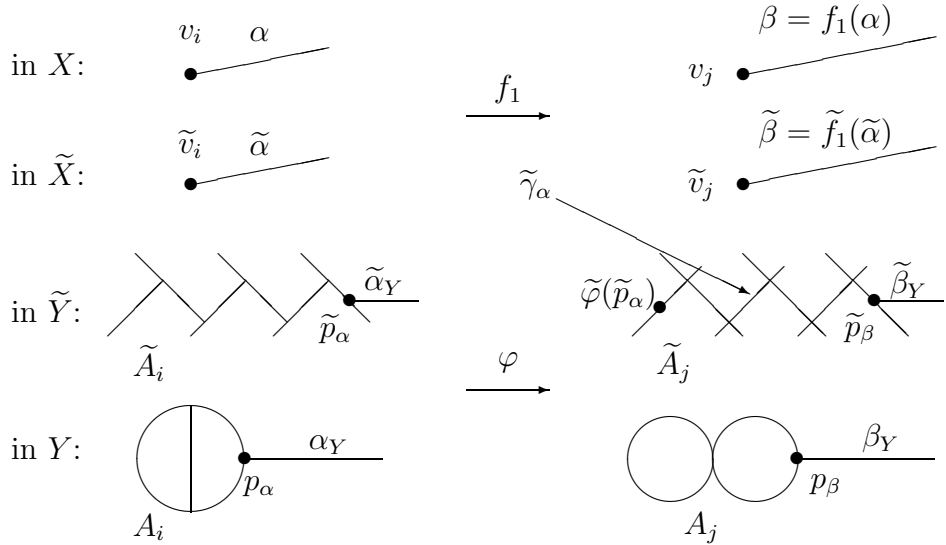

Let $\beta=f_1(\alpha)$ and choose $\wt\beta$ to be the lift of $\beta$ so that $\wt
f_1(\wt\alpha)=\wt\beta$. Note that in case $f_1(\alpha)=\alpha$, $\wt\beta$ is not necessarily
equal to $\wt\alpha$ (it is only in the same orbit). 

Clearly $\wt\beta$ emanates from a lift $\wt v_i$ of $v_i$ so that
$\wt f_1(\wt x)=\wt v_i$.
 The germ $\wt\beta$ corresponds to a germ $\wt\beta_Y$ incident to $\wt A_i$ at a point $\wt
 p_\beta$, where $\wt A_i$ is the component of the preimage of $A_i$ so that $\wt\f(\wt
 A_k)=\wt A_i$. We define $\wt \gamma_\alpha$ as the unique geodesic path in $\wt A_i$ connecting
$\wt\varphi(p_\alpha)$ to $\wt p_\beta$. Now we define $\gamma_\alpha$ as the projection to $Y$
of $\wt\gamma_\alpha$. It is a path from $\f(p_\alpha)$ to $p_\beta$.

\begin{rem}
  We chose a path $\wt \gamma_\alpha$ for any germ $\alpha$ in $X$, which is a finite graph.
Therefore we have only finitely many such $\wt \gamma_\alpha$'s. We can then complete that family of
paths by equivariance. 
\end{rem}

\begin{rem}
  If we use $\wt\alpha$ instead $g\wt\alpha$, then bot $\wt\f(\wt p_\alpha)$ and $\wt p_\beta$
  --- and therefore also $\wt\gamma_\alpha$ --- are translated by $\phi(g)$, hence the path $\gamma_\alpha$ is actually independent on the
  choice of the lift $\wt\alpha$.
\end{rem}

{\bf Case 2}. Let $x\in X$ be such that $f_1(x)=v_i$ for some $i$, but $x$ is not one of the
$v_j$'s. (In case $x$ is not a vertex, up to add $x$ to the simplicial structures of $X$ and
$Y$, so we can consider it as a vertex.) For any
germ of edge  $\alpha$ at $x$ we define $\gamma_\alpha$ as follows. 

First, we fix a base-point  $x_i\in A_i$, and  for any component $\wt A_i$ (of the preimage of
$A_i$) we  choose a lift $\wt x_i\in\wt A_i$. Any germ of edge $\alpha$ at $x$
corresponds to a germ $\alpha_Y$ is $Y$. For any such $\alpha$ we choose a
lift $\wt \alpha$. Since $f_1$ does not collapse edges, $\wt f_1(\wt \alpha)$ is  a germ of
edge $\wt \beta$ at some lift $\wt v_i$ of $v_i$, and
corresponds to a germ $\wt \beta_Y$ at $\wt A_i$ in $\wt Y$. Let $\wt \gamma_\alpha$ be the unique
path in $\wt A_i$ connecting $\wt x_i$ and $\wt \beta_Y$. We finally define $\gamma_\alpha$ as
the projection to $Y$ of $\wt\gamma_\alpha$.

\begin{rem}
As above  we chose only finitely many such $\wt\gamma_\alpha$'s and we can complete the choices
equivariantly.
\end{rem}
\begin{rem}
  The path $\gamma_\alpha$ actually depends on the choices of $x_i$ and $\wt x_i$, but for any
  pair of germs $\alpha_1,\alpha_2$ at $x$, the reduced version of the concatenation
  $\gamma_{\alpha_1}^{-1}\gamma_{\alpha_2}$  does not depend on such choices. 
\end{rem}

Note that, as germs, $\alpha_Y=\sigma_t^{-1}(\alpha)$
and $\beta_Y=\sigma_t^{-1}(\beta)=\sigma_t^{-1}(f_1(\alpha))$.
Now we have a path $\gamma_\alpha\subset A$ for any pre-image of germs at the $v_i$'s, chosen
independently on $t$. Let
$t\in(0,1]$. We define a map $$\overline g_t: Y_t\to Y_t$$ representing $[\phi]$ as follows:
\begin{itemize}
\item in $\sigma_t^{-1}\big(X\setminus f_1^{-1}(\sqcup_i B_\rho(v_i))\big)$ we just set
  $\overline g_t=\sigma_t^{-1}\circ f_1\circ\sigma_t$;
\item in $\sigma_t^{-1}\big(f_1^{-1}(\sqcup_i B_\rho(v_i))\big)\setminus A^t$ we use the paths
  $\gamma_\alpha$. More precisely, let $N$ be a connected component of
  $f_1^{-1}(B_\rho(v_i))$  and let $x\in N$ such that $f_1(x)=v_i$. For any
  edge $\alpha\in N$ emanating from $x$ we define $\overline g_t(\sigma_t^{-1}(\alpha))$ by mapping
  linearly\footnote{I.e. at constant speed} $\sigma_t^{-1}(\alpha)$ to the path given by the
  concatenation of
  $\beta_Y=\sigma_t^{-1}(f_1(\alpha))$ and $\gamma_\alpha$. Note that
  $\overline g_t|_{\sigma_t^{-1}(\alpha)}=\PL(\overline g_t|_{\sigma_t^{-1}(\alpha)})$.
\item in $A^t$ we set $\overline g_t=\varphi$;

\end{itemize}
finally, we set $$g_t=\opt(\PL(\overline g_t))$$
where straightening and optimization are made with respect to the metric structure of $Y_t$.
We now estimate the Lipschitz constant of $\overline g_t$. Clearly we have the lower
bound $$\lambda(Y_t)=\Lip(g_t)\leq \Lip(\overline g_t).$$

Moreover, since on edges of $Y_t\setminus A^t$ the map $\sigma_t$ is just a rescaling of edge-lengths, for any $\e_2>0$ there is $t_{\e_2}>0$ such that $\forall t<t_{\e_2}$
\begin{equation}
  \label{eq:2}
  \Lip(\sigma_t)<1+\e_2 \qquad \Lip(\sigma_t^{-1})<1+\e_2.
\end{equation}

Now we compute an upper bound for $\Lip(\overline g_t)$. 
As $\overline g_t$ is defined in three different regions, namely
\begin{itemize}
\item $\Omega_1=\sigma_t^{-1}\big(X\setminus f_1^{-1}(\sqcup_i B_\rho(v_i))\big)$,
\item $\Omega_2=\sigma_t^{-1}\big(f_1^{-1}(\sqcup_i B_\rho(v_i))\big)\setminus A^t$,
\item $\Omega_3=A^t$;
\end{itemize}
we will estimate $\Lip(\overline g_t)$ on these three regions separately.

In $\Omega_1$ we have
  $\overline g_t=\sigma_t^{-1}\circ f_1\circ\sigma_t$. Then
$$\Lip(\overline g_t|_{\Omega_1})\leq \Lip(\sigma_t^{-1})\Lip(f_1)\Lip(\sigma_t).\footnote{Note that
  $\Lip(\sigma_t)$ and $\Lip (\sigma_t^{-1})$ are not the inverse of each other because
  different edges are stretched by $\sigma_t$ by a priori different amounts.}$$

Hence, by~\eqref{eq:1},~\eqref{eq:2}, and by setting $(1+\e_2)^2(1+\e_1)=1+\e_3$, we have
\begin{equation}
  \label{eq:3}
  \Lip(\overline g_t|_{\Omega_1})\leq(1+\e_2)^2\lambda(X)(1+\e_1)=(1+\e_3)\lambda(X).
\end{equation}
Now, we switch to $\Omega_2$. Let $N$ be a connected component of $f_1^{-1}(\sqcup_i B_\rho(v_i))$. Let $x\in N$ such
that $f_1(x)=v_i$ and let $\alpha$ be an edge of $N$ emanating from $x$.
By definition $\overline g_t$ is linear on $\sigma_t^{-1}(\alpha)$, thus in order to estimate its Lipschitz
constant we need to know only the lengths of $\sigma_t^{-1}(\alpha)$ and  its image.
Clearly $$L_X(\alpha)=L_X(\sigma_t(\sigma_t^{-1}(\alpha)))
\leq \Lip(\sigma_t) L_{Y_t}(\sigma_t^{-1}(\alpha))
\quad\text{ and thus }\quad L_{Y_t}(\sigma_t^{-1}(\alpha))\geq\frac{L_X(\alpha)}{\Lip(\sigma_t)}.$$
Moreover, since we have $L_X(f_1(\alpha))=\rho$, we get
$$\rho\leq \Lip(f_1) L_X(\alpha) \qquad\text{ and so }
\qquad L_X(\alpha)\geq\frac{\rho}{\Lip(f_1)} 
$$
whence, by~\eqref{eq:2} and~\eqref{eq:1}, we obtain
$$L_{Y_t}(\sigma_t^{-1}(\alpha))\geq\frac{\rho}{\Lip(\sigma_t)\Lip(f_1)}
>\frac{\rho}{(1+\e_2)\Lip (f_1)}>
\frac{\rho}{\lambda(X)(1+\e_1)(1+\e_2)}.$$
Since $\gamma_\alpha$ is the same path in $A$ for every $t$, its length in $A^t$ depends
linearly on $t$, namely here is a constant $C_\alpha$ such that $$L_{Y_t}(\gamma_\alpha)= C_\alpha t$$
whence, setting $C=\max_\alpha C_\alpha$,
\begin{eqnarray*}
  &~&\Lip(\overline g_t|_{\sigma_t^{-1}(\alpha)})
\leq\frac{L_{Y_t}(\sigma_t^{-1}(f_1(\alpha))+L_{Y_t}(\gamma_\alpha)}{L_{Y_t}(\sigma_t^{-1}(\alpha))}
\leq \frac{\Lip(\sigma_t^{-1})\rho+tC}{L_{Y_t}(\sigma_t^{-1}(\alpha))}\\
  &<&((1+\e_2)\rho+tC)\frac{\lambda(X)(1+\e_1)(1+\e_2)}{\rho}\\
  &<&(1+\e_2)(\rho+tC)\frac{\lambda(X)(1+\e_1)(1+\e_2)}{\rho}\\
  &=& \lambda(X)(1+\e_3)(1+\frac{tC}{\rho}).
\end{eqnarray*}
Therefore $\forall\e_4>0\exists t_{\e_4}>0$ such that $\forall t<t_{\e_4}$ and for any $\alpha$ as
above, we have $\Lip(\overline g_t|_{\sigma_t^{-1}(\alpha)})<\lambda(X)(1+\e_4)$ and so

\begin{eqnarray}
  \label{eq:4}
  \Lip(\overline g_t|_{\Omega_2})=\sup_\alpha\Lip(\overline g_t|_{\sigma_t^{-1}(\alpha)})<\lambda(X)(1+\e_4).
\end{eqnarray}

Finally,  on $\Omega_3=A^t$ we have $\overline g_t=\varphi$ and so
$\Lip(\overline g_t|_{A^t})=\Lip(\varphi)$. Thus, by~\ref{eq:fi}
\begin{eqnarray}
  \label{eq:7}
\Lip(\overline g_t|_{\Omega_3})\leq \lambda(Y_t).
\end{eqnarray}

Since by~\eqref{eq:0} $\lambda(X)\leq\lambda(Y_t)(1+\e_0)$, by putting
together~\eqref{eq:3},~\eqref{eq:4}, and~\ref{eq:7} we have that for any $\e_5>0$ there is
$t_{\e_5}>0$ such that for any $t<t_{\e_5}$ we have
$$\Lip(\overline g_t)\leq\lambda(Y_t)(1+\e_5).$$
We are now in position to obtain the inequality claimed in the statement.
Since $g_t$ is optimal, $\Lip(g_t)=\lambda(Y_t)$, and by Theorem~\ref{Lemma_opt}
$$d_\infty(g_t,\overline g_t)<\vol(Y_t)(\Lip(\overline g_t)-\lambda(Y_t))<\vol(Y_t)\lambda(Y_t)\e_5.$$

We first estimate $$d_\infty(\sigma_t\circ \overline g_t,f_1\circ \sigma_t).$$
 In $\sigma_t^{-1}\big(X\setminus f_1^{-1}(\sqcup_i B_\rho(v_i))\big)$ we have
  $\overline g_t=\sigma_t^{-1}\circ f_1\circ\sigma_t$ so here the distance is zero. On $A^t$,
  since $\overline g_t(A)=A$,
  for any $i$ there is $j$ such that  we have $\sigma_t(\overline g_t(A_i))=\sigma_t(A_j)=v_j=f_1(v_i)$,
  hence also in $A^t$ the distance is zero.
Finally, let $N$ be a connected component of $f_1^{-1}(\sqcup_i B_\rho(v_i))$. Let $x\in N$ such
that $f_1(x)=v_i$ and let $\alpha$ be an edge of $N$ emanating from
$x$.
The path $\overline g_t(\sigma_t^{-1}(\alpha))$ is given by the concatenation of $\sigma_t^{-1}(f_1(\alpha))$ with
$\gamma_\alpha$. The latter is collapsed by $\sigma_t$, and the image of the former is just
$f_1(\alpha)=f_1\circ\sigma_t(\sigma_t^{-1}(\alpha))$. Since the length of $\gamma_\alpha$ in
$A^t$  is bounded by $tC$ we have that
$$d_\infty(\sigma_t\circ \overline g_t,f_1\circ \sigma_t)\to 0\qquad\text{ as }t\to 0.$$
In particular $\forall \e_6\exists t_{\e_6}$ such that $\forall t<t_{\e_6}$ we have
$$d_\infty(\sigma_t\circ \overline g_t,f_1\circ \sigma_t)<\e_6.$$

Finally,
\begin{eqnarray*}
&~&d_\infty(\sigma_t\circ g_t,f\circ \sigma_t)\\&\leq&
d_\infty(\sigma_t\circ g_t,\sigma_t\circ \overline g_t)+
d_\infty(\sigma_t\circ \overline g_t,f_1\circ \sigma_t)+
d_\infty(f_1\circ\sigma_t,f\circ \sigma_t)\\&\leq&
\Lip(\sigma_t)d_\infty(g_t,\overline g_t)+\e_6+d_\infty(f_1,f)\\&<&(1+\e_2)\vol(Y_t)\lambda(Y_t)\e_5+\e_6+\e_1
\end{eqnarray*}
which is arbitrarily small for $t\to 0$.\qed

\begin{lem}[When $\core(A)$ is trivial]\label{lemmaZ}
Let $X,Y\in\overline{\Og(\Gamma)}$. Suppose that as graphs of groups, $X$ is obtained from $Y$
by collapsing a sub-forest $A=\sqcup A_i$ whose tree $A_i$ each contains at most one non-free
vertex.
Then the conclusion of Theorem~\ref{thmZnew} holds.
\end{lem}
\proof
Except the definition of $g_t$, the proof goes exactly as that of Lemma~\ref{lemmaZnew}, and it
is even simpler. So let's define $g_t$.
As above $A^t$ denote the scaled version of $A$.
Let $v_i$ be the vertex of $X$ resulting from the collapse of $A_i$.
The function $\lambda$ is now continuous
$$\lambda(Y_t)\to \lambda(X).$$
As above, if $f$ collapses some edge we find $f_1:X\to X$ a straight map representing $[\phi]$ which
collapses no edge and with
$$d_\infty(f,f_1)<\e_1\quad \text{and}\quad \Lip(f_1)<\Lip(f)(1+\e_1)=\lambda(X)(1+\e_1).$$

We choose $\rho$ so that $B_\rho(v_i)$ is star-shaped, the components of
$f_1^{-1}(B_\rho(v_i))$ are star-shaped and contain a unique pre-image of $v_i$, and so that
the  components of $f_1^{-1}(B_\rho(v_i))$  and $f_1^{-1}(B_\rho(v_j))$ are pairwise
disjoint. Finally we chose $\rho$ small enough so that if $f(v_i)\notin\{v_j\}$,
then $f(v_i)\notin\cup_jB_\rho(v_j)$.

For any $i$ we choose a base vertex $x_i\in A_i$ which is the non-free vertex of $A_i$ if
any. For any $x\in X$ such that $f_1(x)=v_i$ and for any edge $\alpha$ in $f_1^{-1}(B_\rho(v_i))$
incident to $x$, let $\gamma_\alpha$ be the unique embedded path connecting
$\sigma_t^{-1}(f_1(\alpha))$ to $x_i$. We define $\overline g_t:Y_t\to Y_t$ as follows:
\begin{itemize}
\item in $\sigma_t^{-1}\big(X\setminus f_1^{-1}(\sqcup_i B_\rho(v_i))\big)$ we just set
  $\overline g_t=\sigma_t^{-1}\circ f_1\circ\sigma_t$;
\item in $\sigma_t^{-1}\big(f_1^{-1}(\sqcup_i B_\rho(v_i))\big)\setminus A^t$ we use the paths
  $\gamma_\alpha$. More precisely, let $N$ be a connected component of
  $f_1^{-1}(B_\rho(v_i))$  and let $x\in N$ such that $f_1(x)=v_i$. For any
  edge $\alpha\in N$ emanating from $x$ we define $\overline g_t(\sigma_t^{-1}(\alpha))$ by mapping
  linearly\footnote{I.e. at constant speed} $\sigma_t^{-1}(\alpha)$ to the path given by the
  concatenation of
  $\sigma_t^{-1}(f_1(\alpha))$ and $\gamma_\alpha$. Note that
  $\overline g_t|_{\sigma_t^{-1}(\alpha)}=\PL(\overline g_t|_{\sigma_t^{-1}(\alpha)})$.
\item in the components $A_i^t$ so that $f_1(v_i)=v_j$, we set $g(A_i^t)=x_j$;
\end{itemize}
finally we set $g_t=\opt(\PL(\overline g_t))$. The estimates on Lipschitz constants and distances now
follow exactly as in the proof of Lemma~\ref{lemmaZnew}.\qed

\providecommand{\bysame}{\leavevmode\hbox to3em{\hrulefill}\thinspace}
\providecommand{\MR}{\relax\ifhmode\unskip\space\fi MR }
\providecommand{\MRhref}[2]{%
  \href{http://www.ams.org/mathscinet-getitem?mr=#1}{#2}
}
\providecommand{\href}[2]{#2}

\end{document}